\documentclass[]{interact}

\usepackage{epstopdf}
\usepackage[caption=false]{subfig}

\usepackage[numbers,sort&compress]{natbib}
\usepackage{verbatim}
\usepackage{hyperref}
\usepackage{booktabs}
\usepackage{siunitx}
\usepackage{fancyvrb}
\usepackage{graphicx}
\usepackage{comment}
\usepackage{mathtools}
\usepackage{listings}
\usepackage{xcolor}
\usepackage[ruled,vlined]{algorithm2e}

\DeclareMathOperator*{\st}{subject\,to}
\newcommand{\norm}[1]{\left\lVert#1\right\rVert}
\newcommand{\xx}{\mathbf{x}}
\newcommand{\yy}{\mathbf{y}}

\newcommand{\Ab}{\mathbf{A}}

\newcommand{\X}{\mathbf{X}}
\newcommand{\bb}{\mathbf{b}}

\newcommand{\thh}{\pmb{\theta}}
\newcommand{\Gam}{\pmb{\Gamma}}

\newcommand{\RR}{\mathbb{R}}
\newcommand{\absv}[1]{\vert#1\vert}
\newcommand{\vb}[1]{{\Verb!#1!}}

\bibpunct[, ]{[}{]}{,}{n}{,}{,}
\renewcommand\bibfont{\fontsize{10}{12}\selectfont}

\theoremstyle{plain}

\theoremstyle{definition}

\theoremstyle{remark}
\newtheorem{remark}{Remark}

\definecolor{codegreen}{rgb}{0,0.6,0}
\definecolor{codegray}{rgb}{0.5,0.5,0.5}
\definecolor{codepurple}{rgb}{0.58,0,0.82}
\definecolor{backcolour}{rgb}{0.99,0.99,0.99}

\lstdefinestyle{mystyle}{
    backgroundcolor=\color{backcolour},   
    commentstyle=\color{codegreen},
    keywordstyle=\color{red},
    numberstyle=\color{blue},
    stringstyle=\color{codepurple},
    basicstyle=\ttfamily\footnotesize,
    breakatwhitespace=false,         
    breaklines=true,                 
    captionpos=b,                    
    keepspaces=true,                 
    showspaces=false,                
    showstringspaces=false,
    showtabs=false,                  
    tabsize=2
}

\begin{document}

\articletype{Article}

\title{Sparse Convex Optimization Toolkit: A Mixed-Integer Framework}

\author{
\name{Alireza Olama\textsuperscript{a}, Eduardo Camponogara\textsuperscript{a}, and Jan Kronqvist\textsuperscript{b}}
\affil{\textsuperscript{a}Automation and System Engineering Department, Federal University of Santa Catarina, Florian\'opolis, Brazil\\\textsuperscript{b}Department of Mathematics, KTH Royal Institute of Technology, Stockholm, Sweden}
}

\maketitle

\begin{abstract}
This paper proposes an open-source distributed solver for solving Sparse Convex Optimization (SCO) problems over computational networks. 
Motivated by past algorithmic advances in mixed-integer optimization, the Sparse Convex Optimization Toolkit (SCOT) adopts a mixed-integer approach to find exact solutions to SCO problems. 
In particular, SCOT brings together various techniques to transform the original SCO problem into an equivalent convex Mixed-Integer Nonlinear Programming (MINLP) problem that can benefit from high-performance and parallel computing platforms.
To solve the equivalent mixed-integer problem, we present the Distributed Hybrid Outer Approximation (DiHOA) algorithm that builds upon the LP/NLP based branch-and-bound and is tailored for this specific problem structure. The DiHOA algorithm combines the so-called single- and multi-tree outer approximation, naturally integrates a decentralized algorithm for distributed convex nonlinear subproblems, and utilizes enhancement techniques such as quadratic cuts. Finally, we present detailed computational experiments that show the benefit of our solver through numerical benchmarks on 140 SCO problems with distributed datasets. To show the overall efficiency of SCOT we also provide performance profiles comparing SCOT to other state-of-the-art MINLP solvers.
\end{abstract}


\begin{table}
\tbl{Acronyms}
{\begin{tabular}{ll} \toprule
  Acronym & Meaning   \\ \midrule
 API &  Application Programming Interface  \\
 BnB  & Branch and Bound \\
 CLI & Command Line Interface \\
 DiHOA & Distributed Hybrid Outer Approximation\\
 DiPOA & Distributed Primal Outer Approximation\\
 D-MINLP & Distributed Mixed Integer Nonlinear Programming\\
 DSLinR &   Distributed Sparse Linear Regression \\
 DSLogR & Distributed Sparse Logistic Regression  \\
 LFC  & Local Fusion Center \\
 MILP & Mixed Integer Linear Programming\\
 MINLP & Mixed Integer Nonlinear Programming\\
 MIP & Mixed Integer Programming\\
 MPI & Message Passing Interface \\
 NLP & Nonlinear Programming\\
 OA & Outer Approximation\\
 QCLP & Quadratically Constrained Linear Programming\\
 RH-ADMM & Relaxed Hybrid Alternating Direction Method of Multipliers\\
 SCO & Sparse Convex Optimization\\
 SCOT & Sparse Convex Optimization Toolkit \\
 SOS-1 & Specially Ordered Set of Type I  \\ \bottomrule
\end{tabular}}
\label{}
\end{table}

\section{Introduction}\label{sec:introduction}
In recent years, \textit{Sparse Convex Optimization (SCO)} has gained considerable attention in several disciplines, from machine learning and engineering to economics and finance \cite{tillmann2021cardinality,bertsimas2021sparse,Bertsimas2017a}.  
Several mathematical optimization problems in this context can be formulated as a general convex optimization problem subject to a constraint that allows only up to a certain number of decision variables to be nonzero.
We refer to this constraint as a \textit{sparsity constraint}. Hence, any convex optimization problem with the sparsity constraint can be regarded as a SCO problem \cite{bienstock1996computational,bertsimas2021sparse,bertsimas2022scalable,Sun2013}.

Due to the non-convexity and discontinuity of the sparsity constraint, the SCO problems are known to be NP-Hard \cite{natarajan1995sparse,bertsimas2021sparse}. To overcome the computational difficulties imposed by the sparsity constraint,  computationally tractable convex optimization-based methods have been proposed. One of the popular methods is a $\ell_1$ norm relaxation method where a $\ell_1$ regularizer is imposed on the decision vector. The $\ell_1$ method naturally produces a sparse solution by setting many variables to zero. One of the popular $\ell_1$ based methods is Lasso which is widely used in statistics and machine learning communities \cite{boyd2011distributed,tibshirani1996regression,chen2001atomic}.

One of the important reasons behind the popularity of $\ell_1$ based methods is their computational efficiency and scalability to practical-sized problems. However, in spite of their favorable computational properties, these methods can have some shortcomings. For example, they cannot guarantee that $\ell_1$ based methods find the correct sparsity for general problems. Moreover, in some applications, the desired sparsity structure is different from general sparsity and cannot be easily obtained by a $\ell_1$ regularization. An example of such a sparse structure is \textit{group sparsity} in which a block or a group of independent variables are either all zero or all nonzero. Some notable applications with group sparsity are block-wise linear regression \cite{kim2006blockwise}, logistic regression  \cite{Bertsimas2017a}, compressed sensing \cite{eldar2012compressed}, and microarray analysis \cite{ma2007supervised}. 

Another approach to solving the SCO problems is to view the SCO problems as equivalent Mixed Integer Programming (MIP) problems. Considering recent advances in mixed-integer optimization algorithms and technologies, the MIP problems can be solved efficiently by current mixed-integer optimization solvers such as Gurobi \cite{gurobi}. The resulting MIP formulation is flexible and can be adjusted based on the application's needs. Moreover, by defining suitable binary variables, the MIP framework can provide  exact sparse solutions to SCO problems. 
The MIP framework is gaining popularity in various areas such as statistical data analysis and interpretable machine learning \cite{bertsimas2021sparseclass,bertsimas2016best,bertsimas2021sparse2,bertsimas2021sparse}, sparse control \cite{aguilera2017quadratic}, unit commitment, face recognition, sensor network design \cite{lewis2004wireless}, portfolio optimization \cite{bertsimas2022scalable}, and compressed sensing \cite{foucart2013invitation}. 

The mentioned works focused on \textit{centralized} solutions to SCO problems that might not be suitable for modern real-world applications when the data is inherently distributed or available in large volumes.
Considering the limitations of centralized architectures, in the past few years, mainly because of the rise of Big data, distributed optimization over networks has gained growing attention \cite{notarstefano2019distributed,notarnicola_distributed_2017}. The primary purpose of distributed optimization is to solve an optimization problem over a network of computing nodes. 
Each node performs local computations with access only to a portion of the problem data. Nodes are also capable of exchanging information with other nodes in the network. See \cite{notarstefano2019distributed} for a comprehensive overview of the most common distributed optimization algorithms. 
%

This paper introduces \vb{SCOT}, a distributed optimization solver designed to solve SCO problems over peer-to-peer networks of computing nodes. In essence, \vb{SCOT} consists of two distributed algorithms developed by the authors to solve SCO problems where an iterative procedure is applied by each network node, alternating communication, and computation phases until a solution is found. 
To the best of our knowledge, \vb{SCOT} is the first software framework that can solve SCO problems using distributed algorithms, enabling practical applications with a large number of sample data points, while keeping the data private to each node and allowing implementation in a computer cluster.

Formally, we consider the SCO problem as a mathematical programming problem that consists of finding the $\kappa$-sparse optimal solution of a distributed convex optimization problem of the following form,
\begin{subequations}\label{dccp}
\begin{align}
    \min_{\xx \in \RR^n}\,\, &\sum_{i = 1}^N f_i(\xx)\\
    \st~ & \Ab \xx \leq \bb\\
         & \norm{\xx}_0 \leq \kappa \label{card-constr} 
\end{align}
\end{subequations}
where $N$ is the number of nodes of the computation network, $\xx \in \RR^n$ is the vector of decision variables, and $f_i: \RR^n \xrightarrow{} \RR$ is a convex function assumed to be continuously differentiable and only known by node $i$, for all $i \in \{1,\dots, N\}$. We use the  $\ell_0$ norm (\textit{i.e.}, $\norm{\xx}_0 = \absv{supp(\xx)}= \absv{\{j:x_j \neq 0\}} $) to define the sparsity constraint, which imposes the number of non-zero elements of $\xx$ to be less than a given integer $\kappa$. Finally, the matrix $\Ab \in \RR^{m \times n}$ and the vector $\bb \in \RR^{m}$ define the set of given linear constraints assumed to be known by all nodes.

\subsection{Main Contributions}
The main idea behind \vb{SCOT} is to transform problem \eqref{dccp} into an equivalent \textit{Distributed Mixed Integer Nonlinear Programming (D-MINLP)} problem.  By using the concept of the so-called Local Fusion Centers (LFCs) presented in  \cite{olama2019}, we recast problem \eqref{dccp} as a consensus optimization problem and introduce several constraints to model the sparsity constraint \eqref{card-constr}.
Based on the Outer Approximation (OA) algorithm \cite{Kronqvist2020, Grossmann2002, duran1986outer}, \vb{SCOT} consists of two main algorithms, namely, \textit{Distributed Primal Outer Approximation (DiPOA)} and \textit{Distributed Hybrid Outer Approximation (DiHOA)}. 

The DiPOA algorithm proposed in \cite{olama2021dipoa} extends the OA algorithm by embedding a fully decentralized algorithm, namely the Relaxed Hybrid Alternating Direction Method of Multipliers (RH-ADMM) \cite{olama2019}. The RH-ADMM assumes a particular hybrid architecture on the computational network, developed to solve distributed convex optimization problems. In particular, DiPOA solves the primal problem  of the OA using the RH-ADMM algorithm and handles distributedly the demanding computational part of the OA algorithm that deals with the solution of convex NLPs. In practice, there exist many cases where a significant part of the solution time is spent on solving the NLP problems.
For example, in Table 1 \cite{Kronqvist2020}, it  can  be seen that OA spends more than $150$ seconds on solving NLP problems when solving a moderate-size convex MINLP. Moreover, for inherently distributed problems in which the data is spread over a possibly large computational network, a single NLP will be challenging, if even possible, to solve in a classical centralized fashion. For such inherently decentralized problems, a distributed algorithm can offer great computational advantages. 

Despite solving problem \eqref{dccp} distributedly,  DiPOA is developed based on a multiple-tree OA algorithm whereby a BnB tree is built from scratch at each iteration of the algorithm. Therefore, most of the overall solution time of DiPOA is usually spent on solving MIP sub-problems. In this paper, we tackle the limitations of DiPOA by proposing DiHOA, a distributed algorithm that improves DiPOA performance by gradually building a single BnB tree to avoid constructing and solving many similar MILP problems from scratch.
The main idea of DiHOA is to solve problem \eqref{dccp} by dynamically updating the MILP subproblem, in a fashion similar to the LP/NLP-BnB  presented by \cite{quesada1992lp}. In principle, DiHOA starts with the multiple-tree search strategy up to a certain iteration and introduces high-quality cuts until a suitable event is triggered. 
Once the event is triggered, DiHOA switches from the multiple-tree search strategy to the single-tree search strategy that builds a single BnB tree, starting with multiple cuts initially introduced to its root node to augment the initial formulation. The single BnB tree then tightens up the integer relaxations by dynamically introducing more linear approximations (cuts) to the MILP problem. The multiple-tree search strategy is only applied in some initial iterations of the DiHOA since the MIP problems are typically much easier to solve, as the nonlinear constraints are only roughly represented through a few constraints. Moreover starting the BnB search with a tighter approximation of the nonlinear constraints can result in a much smaller BnB tree, as a smaller infeasible region of the continuously relaxed search space is explored. Without these initial cuts, the nonlinear constraints would be completely ignored until an integer solution is found and the nonlinear constraints would be poorly approximated until a few integer solutions is explored in the BnB tree.

The main contributions of this work are summarized as follows:
\begin{itemize}

    \item \vb{SCOT}, a distributed software framework to model and solve SCO problems.
    
    \item The distributed hybrid outer approximation algorithm to solve the SCO problem \eqref{dccp}.
    
    \item Modeling and heuristic techniques to improve the efficiency of the proposed algorithm.
    
    \item A computational analysis of the proposed algorithms for various SCO real-world applications.
\end{itemize}
\subsection{Paper Organization}
The paper is organized as follows. In Section \ref{sec:problem} we introduce the distributed SCO problem by reformulating problem \eqref{dccp} into an equivalent MINLP problem. Section \ref{sec:primaldual} presents the primal and dual problems used by \vb{SCOT} algorithms. The DiHOA algorithm is proposed in Section \ref{sec:algorithms}. Section \ref{sec:scot} introduces \vb{SCOT} and its main components. Finally, the numerical comparison and algorithm analysis are presented in Section \ref{sec:applications}. Section \ref{sec:conclusion} concludes the paper.

\section{Distributed Sparse Convex Optimization}\label{sec:problem}
Various modeling techniques are used by \vb{SCOT} to find a solution to problem \eqref{dccp} efficiently. First, \vb{SCOT} transforms problem \eqref{dccp} into a consensus optimization problem and then multiple modeling techniques are used to handle the sparsity constraint \eqref{card-constr}.

\subsection{Consensus Optimization Modeling}
As defined in problem \eqref{dccp}, the decision variables are shared between the nodes. Each node only has information to construct its objective function, while keeping the local problem data private from other nodes.
This distributed setting is a typical pattern in some significant learning and control applications. For instance, in distributed machine learning, an efficient technique to deal with large volumes of data consists of distributing the data over a network \cite{notarstefano2019distributed,nedic2018distributed}. 
In this case, a significant reduction in memory size for computation can be achieved while keeping the same unknown parameters of the model. 

To decompose \eqref{dccp}, we adhere to the concept of \textit{hypergraphs}, which is a generalization of a regular graph in which an edge can join an arbitrary number of vertices. Multiple computational sources, called LFCs, can be employed by adopting this structure in the network. 
We consider a hypergraph $\mathcal{H}=(\mathcal{V},\mathcal{E})$ defined as follows: $\mathcal{V}=\{1,2, ...,N\}$ is the set of nodes such that node $i$ decides upon the values of vector variable $\mathbf{x}_i$; $\mathcal{E}=\{\mathcal{E}_k\subset \mathcal{V} : k=1,\ldots,K\}$ is the set of hyperedges, where a hyperedge  $\mathcal{E}_k$ connects all nodes $i\in\mathcal{E}_k$ and $K$ is the number of hyperedges. Now we introduce the concept of path in a hypergraph: a path $p(i,j)=\langle \mathcal{E}_1,..., \mathcal{E}'_k\rangle$ connects nodes $i$ and $j$ if $i\in \mathcal{E}'_1$, $j\in \mathcal{E}'_k$, and $\mathcal{E}'_l\cap \mathcal{E}'_{l+1}\neq \emptyset$ for $l=1,\ldots,k-1$, and $\mathcal{E}'_l\in\mathcal{E}$ for all $l$. 
Put another way, through the hyperedges, a path $p(i,j)$ establishes a communication channel between nodes $i$ and $j$.

We assume that the hypergraph $\mathcal{H}$ is connected, meaning that for all $i,j\in\mathcal{V}$ there exists a path $p(i,j)$ connecting $i$ and $j$. By using the hypergraph structure, the equivalent formulation for \eqref{dccp} is obtained as follows,
\begin{subequations}\label{dcco-lfc}
\begin{align}
    \min_{\substack{\xx_1,...,\xx_N \\ \yy_1,...,\yy_K}}\,\, &\sum_{i = 1}^N f_i(\xx_i)\\
    \st~ & \Ab \xx_i \leq \bb,  \, \forall i=1,...,N,\\
         & \xx_i = \yy_j, \, \forall i \in \mathcal{E}_j, \, \mathcal{E}_j \in \mathcal{E} \\
         & \norm{\yy_j}_0 \leq \kappa, \, \forall j=1,\ldots,K \label{card-constr-lfc} 
\end{align}
\end{subequations}
\textcolor{black}{where $\xx_i \in \RR^n$ are vectors of decision variables associated with the nodes and $\yy_j \in \RR^n$ are auxiliary variables associated with the LFCs, which are represented by the hyperedges.}

\subsection{Sparsity Constraint Modeling}
Here, we present three modeling techniques implemented by \vb{SCOT} to express the sparsity constraints \eqref{card-constr-lfc}, namely the \textit{big-M}, \textit{Specially Ordered Set of Type I} (SOS-1), and a \textit{hybrid} approach. 

\subsubsection{Big-M Method}
The Big-M method is arguably the simplest technique for modeling the sparsity constraint. This method incorporates a binary variable and an estimated upper bound into the model for each continuous variable appearing in the sparsity constraint. By using the Big-M method, we recast the sparsity constraint \eqref{card-constr-lfc} as the following inequalities,
\begin{subequations}\label{bim-m}
\begin{align}
         & -{\tt M}_j\delta_{jk} \leq y_{jk} \leq {\tt M}_j \delta_{jk}, \, \forall j=1,\dots,K, \, \forall k=1,\dots,n \label{bigm} \\
         &\textcolor{black}{\sum_{k=1}^{n} \delta_{jk} \leq \kappa, \, \forall j=1,\dots,K} \label{sos1}\\
         & \delta_{jk} \in \{0, 1\}, \, \forall j=1,\ldots,K, \, \forall k=1,\dots,n 
\end{align}
\end{subequations}
where $y_{jk}$ is the $k$-th element of $\yy_j$, $\boldsymbol{\delta}_{j}, \, \forall j=1,\ldots,K$, is a vector of binary variables whose $k$-th element is denoted by $\delta_{jk}$, and ${\tt M}_j$ is a constant assumed to be a valid upper bound for $\norm{\yy_j}_{\infty}$. In this case, if $\delta_{jk} = 0$ then  $y_{jk}=0$ and $y_{jk}=1$ otherwise. Thus, inequalities \eqref{bigm} and \eqref{sos1} impose the maximum number of nonzero variables in $\yy_j$.
 The big-M parameter ${\tt M}_j$ is not known a priori, and a too small value of ${\tt M}_j$ may lead to a sub-optimal solution. A large ${\tt M}_j$ on the other hand, will result in a weak continuous relaxation and strongly affect the number of nodes that need to be explored in BnB. Hence, a good choice of ${\tt M}_j$ affects the strength of the formulation, being critical for MIP algorithms to obtain high-quality bounds. In the context of learning and control applications, the big-M value ${\tt M}_j$ can typically be computed from data in statistical learning tasks \cite{bertsimas2021sparseclass} and from the physical bounds in control applications \cite{aguilera2017quadratic}. 

\subsubsection{Specially Ordered Set of Type I (SOS-1) Method}
This section discusses the sparsity constraint reformulation using the SOS-1 constraint. Any feasible solution to problem \eqref{dcco-lfc} satisfies the following complementary constraints,
\begin{equation}\label{pre-sos}
    (1 - \delta_{jk})y_{jk} = 0, \, \forall j=1,\ldots,K, \, \forall k=1,\dots,n
\end{equation}
\textcolor{black}{which is equivalent to constraint \eqref{card-constr-lfc}.}
In order for  constraint \eqref{pre-sos} to be satisfied, either $(1 - \delta_{jk})$ or $y_{jk}$ must be zero. Such constraints can be modeled via integer optimization software using Specially Ordered Sets of Type I (SOS-1) \cite{bertsimas2005optimization} as follows,
\begin{equation}
    \left ( y_{jk}, 1 - \delta_{jk}\right )\,:\,\text{SOS-1}, \, \forall j=1,\ldots,K, \, \forall k=1,\dots,n.
\end{equation}
It is worth noting that such a constraint is not implemented explicitly employing algebraic equations, as in the Big-M method. Instead, the optimization solver can directly enforces the SOS-1 constraint by branching on sets of variables. The MIP solvers usually adopt different approaches to handle the SOS-1 constraints. One approach is to use an equivalent big-M formulation in case strong bounds on variables can be obtained. Otherwise, the MIP solver can also handle the constraint directly through branching. By using the SOS-1 technique, we allow the MIP solver to choose the approach it determines to be the best choice to deal with the SOS-1 constraints.

\subsubsection{MINLP Reformulation}
This section discusses the primary problem formulation that \vb{SCOT} considers. Although the big-M formulation is relatively straightforward, choosing a suitable big-M value ${\tt M}_j$ is important. In case ${\tt M}_j$ is too small, the optimization algorithm may cut off valid solutions. However, if ${\tt M}_j$ is excessively large, the model may become numerically difficult to solve. SOS-1 constraints, on the other hand, have the advantage of avoiding these types of problems, as they do not rely on a problem-dependent constant value. For applications that require a relatively small value of ${\tt M}_j$, the big-M modeling technique may be a proper choice, although, for other cases, the SOS-1 method may be suitable. In this regard, \vb{SCOT} can provide the functionality to choose the most proper method depending on the application and problem data. Considering big-M and SOS-1 constraints for modeling the sparsity constraint, the main optimization problem that \vb{SCOT} attempts to solve is expressed as follows,
\begin{subequations}\label{dcco-lfc-sos1-bigm-epig}
\begin{align}
    \min_{\substack{\gamma \\ \xx_1,...,\xx_N \\ \yy_1,...,\yy_K } }\,\, &\gamma \\
    \st~ &\sum_{i = 1}^N f_i(\xx_i) - \gamma\leq 0 \\
         &\Ab \xx_i \leq \bb,  \, \forall i=1,...,N,\\
         & \xx_i = \yy_j, \, \forall i \in \mathcal{E}_j, \, \mathcal{E}_j \in \mathcal{E} \\  
         & -{\tt M}_j\textcolor{black}{\delta_{jk}}\leq y_{jk} \leq {\tt M}_j\textcolor{black}{\delta_{jk}}, \, \forall j=1,\ldots,K, \, \forall k=1,\ldots,n \label{bigm-minlp}\\  
         & \left(y_{jk}, 1 - \delta_{jk}\right )\,:\,\text{SOS-1}, \, \forall j=1,\ldots,K, \, \forall k=1,\ldots,n \label{sos1-minlp}\\  
         &\textcolor{black}{\sum_{k=1}^{n} \delta_{jk} \leq \kappa, \, \forall j=1,\dots,K} \label{mip-const-sos}\\
         & \delta_{jk} \in \{0, 1\}, \, \forall j=1,\ldots,K, \, \forall k=1,\ldots,n 
\end{align}
\end{subequations}
in which the equivalent epigraph reformulation is used. The advantage of using both constraints \eqref{bigm-minlp} and \eqref{sos1-minlp} is that it might be possible to determine better ${\tt M}_j$ coefficients than the MIP solver is able to derive. In that way, \vb{SCOT} provides more information to the MIP solver by including constraint \eqref{bigm-minlp} and \eqref{sos1-minlp}. Problem \eqref{dcco-lfc-sos1-bigm-epig} is an MINLP problem with a separable structure, where the objective function is linear and $\gamma \in \mathcal{R}$ is an auxiliary variable. In the following sections, we introduce a distributed formulation and algorithm that \vb{SCOT} implements to solve problem \eqref{dcco-lfc-sos1-bigm-epig}.

\section{SCOT Primal and Dual Problem}\label{sec:primaldual}
Here, we present the dual and primal problems and we discuss two algorithms that \vb{SCOT} implements in the next section. Akin to other decomposition-based MINLP algorithms, \vb{SCOT} decomposes problem \eqref{dcco-lfc-sos1-bigm-epig} into two main sub-problems, namely, the \textit{primal} and \textit{dual} problems \cite{Kronqvist2020, Lundell2019}. We use the term \textit{primal solution} and \textit{primal bound} as the optimal solution and objective value of the primal problem respectively. Similarly, \textit{dual solution} and \textit{dual bound} are used for the dual problem. It should be noted that we use the terminology dual problem for a problem whose optimal solution provides a valid lower bound on the optimal objective value of problem \eqref{dcco-lfc-sos1-bigm-epig} and whose feasible set contains all feasible solutions of problem \eqref{dcco-lfc-sos1-bigm-epig}. 

The primal solution is assumed to satisfy all linear, nonlinear, and consensus constraints of problem \eqref{dcco-lfc-sos1-bigm-epig} to a given tolerance. Additionally, the current best-known primal solution found by the algorithm is referred to as the \textit{incumbent solution}, and its objective value is the current primal bound. The incumbent solution is updated when \vb{SCOT} algorithms find a primal solution with a lower objective value. In the standard OA algorithm, the primal problem is a convex NLP problem; however, owing to the distributed nature of problem \eqref{dcco-lfc-sos1-bigm-epig}, the primal problem of \vb{SCOT} is a distributed convex NLP problem which will be discussed later.

A solution point whose objective value provides a valid lower bound for the optimum of problem \eqref{dccp}, but not necessarily satisfying all constraints, is referred to as a dual solution. Like the standard OA algorithm, \vb{SCOT} obtains dual solutions by solving relaxed problems that approximate the nonlinear constraints with polyhedral outer approximations. Depending on the type of outer approximations, the dual problem can be MILP, MIQP, or Mixed Integer Quadratically Constrained Linear (Quadratic) (MIQCL(Q)P) problems. Moreover, the dual bound is the best possible objective value of the dual problem. 
The primal and dual sub-problems are then iteratively solved by proper MIP algorithms. The MIP algorithms are distinguished depending on how the sub-problems are constructed, solved, and coordinated. Regardless of the solution algorithms adopted by \vb{SCOT}, the dual and primal sub-problems are two primary components of the algorithms.

The main problem reformulation that \vb{SCOT} attempts to solve by default is problem \eqref{dcco-lfc-sos1-bigm-epig}, which enforces both big-M and SOS-1 constraints and epigraph-reformulation. The separability of nonlinear functions in problem \eqref{dcco-lfc-sos1-bigm-epig} allows \vb{SCOT} to employ an alternative formulation, the so-called \textit{lifted formulation} \cite{kronqvist2018reformulations}. In this context, we use the lifted formulation for each nonlinear function, $f_i$, which results in tighter outer approximations when approximating nonlinear functions \cite{kronqvist2018reformulations,hijazi2014outer,tawarmalani2005polyhedral}.

According to the idea of lifted formulation and also following the procedure presented in \cite{kronqvist2018reformulations}, at each iteration $q$, we construct the dual problem as the following MIP problem,
\begin{subequations}\label{dcco-lfc-sos1-bigm-master}
\begin{align}
    \min_{\substack{\xx_1,...,\xx_N \\ \yy_1,...,\yy_K \\ \gamma_1,...,\gamma_N} }\,\, &\sum_{i = 1}^N \gamma_i\\
    \st~ &  \Tilde{f_i}(\xx_i^q) - \gamma_i \leq  0, \,\, \forall \textcolor{black}{\xx^q_i} \in \mathcal{X}_i^q, \, \forall i=1,...,N, \label{approx}\\
         &\Ab \xx_i \leq \bb,  \, \forall i=1,...,N,\\
         & \xx_i = \yy_j, \, \forall i \in \mathcal{E}_j, \, \mathcal{E}_j \in \mathcal{E} \\
         & -{\tt M}_j\textcolor{black}{\delta_{jk}} \leq y_{jk} \leq {\tt M}_j\textcolor{black}{\delta_{jk}} , \, \forall j=1,\ldots,K, \, \forall k=1,\ldots,n\\
         & \left (y_{jk}, 1 - \delta_{jk}\right )\,:\,\text{SOS-1}, \, \forall j=1,\ldots,K, \, \forall k=1,\ldots,n\\
         &\textcolor{black}{\sum_{k=1}^{n} \delta_{jk} \leq \kappa, \, \forall j=1,\dots,K} \\
         & \delta_{jk} \in \{0, 1\}, \, \forall j=1,\ldots,K, \, \forall k=1,\ldots,n 
\end{align}
\end{subequations}
where $\gamma_i \in \mathcal{R}, \, i = \{1,...,N\},$ are new auxiliary decision variables, $\xx_i^q$ is a feasible point that satisfies linear and consensus constraints, obtained at iteration $q$, and $\Tilde{f_i} (\xx_i^q)$ is the outer approximation of $f_i(\xx_i)$ around $\xx_i^q$.  Finally, $\mathcal{X}_i^q$ is a finite set consisting of local feasible points defined as,
\begin{align*}   
   \mathcal{X}_i^q &= \left \{\xx_i^\ell: \Ab \xx_i^\ell \leq \bb, \, \forall \ell \in \{1,\ldots,q\}\right \},
\end{align*}
Set $\mathcal{X}_i^q$ consists of the feasible points up to the current iteration $q$. The dual problem \eqref{dcco-lfc-sos1-bigm-master} is a relaxation of problem \eqref{dcco-lfc-sos1-bigm-epig} since approximations of nonlinear constraints are used and its objective value is a lower bound of \eqref{dcco-lfc-sos1-bigm-epig}.  At each iteration of \vb{SCOT} algorithms, $q$, a new outer approximation is generated by each node of the network and cooperatively added to problem \eqref{dcco-lfc-sos1-bigm-master} producing a tighter representation of the nonlinear constraints. 
The quality of outer approximations generated by $\Tilde{f_i} (\xx_i^q)$ directly impacts the convergence of the algorithms. Hence, \vb{SCOT} provides first and second-order outer approximations and an event-triggered scheme that controls the effectiveness of the approximations. According to first-order Taylor series and convexity of $f_i(\xx_i)$ functions, we can express constraints \eqref{approx} as,
\begin{equation} \label{foc}
    f_i(\xx_i^q)  + \nabla f_i(\xx_i^q)^T (\xx_i - \xx_i^q) - \gamma_i \leq 0,
\end{equation}
which are linear inequalities. In case the nonlinear functions, $f_i(\xx_i)$, are \textit{strongly} convex functions, \vb{SCOT} replaces constraints \eqref{approx} with the following quadratic inequalities,
\begin{equation} \label{soc}
    f_i(\xx_i^q)  + \nabla f_i(\xx_i^q)^T (\xx_i - \xx_i^q) + \frac{m_i^q}{2}\norm{(\xx_i - \xx_i^q)}_2^2 - \gamma_i \leq 0 
\end{equation}
where $m_i^q > 0 $ is a constant such that $\nabla^2f_i(\xx_i) \succeq m_i^q I$. With $m_i^q>0$, it is clear that the cut given by \eqref{soc} is stronger than the cut given by \eqref{foc}. However, the quadratic cuts \eqref{soc} tend to result in more challenging sub problems in BnB. Therefore, there can still be a computational advantage of the linear cuts.

\begin{remark}
    For general strongly convex functions, $m_i^q$ is not obtained easily. However, in some practical problems found in statistical learning and control, the computation of $m_i^q$ is feasible. For example, the objective function in sparse Model Predictive Control (s-MPC) problems (which is a subclass of the SCO problem) is typically a convex quadratic function. For convex quadratic functions, $m_i^q$ is the smallest Eigenvalue of the Hessian matrix. In machine learning problems the objective function usually consists of a convex function and a strongly convex \textit{regularization} term. In this case, $m_i^q$  can be computed from the regularization term.
\end{remark}
A crucial step to forming the outer approximations is the computation of the approximation points $x_i^q$ for which various strategies and methods exist.  One of the well-known methods to obtain $x_i^q$ is fixing the local binary decision variables of problem \eqref{dcco-lfc-sos1-bigm-epig}, $\delta_{jk} = \delta_{jk}^q$,  and solving the resulting nonlinear optimization problem. In this case, a distributed convex NLP problem is solved, and its optimal solution provides a primal solution to the original MINLP problem. 
The problem of solving \eqref{dcco-lfc-sos1-bigm-epig} for fixed binary variables is the primal problem of \vb{SCOT}, which, at each iteration $q$, is defined as,
\begin{subequations}\label{dcco-lfc-sos1-bigm-primal}
\begin{align}
    \min_{\substack{\gamma \\ \xx_1,...,\xx_N \\ \yy_1,...,\yy_K } }\,\, &\sum_{i=1}^{N}\gamma_i \\
    \st~ & f_i(\xx_i) - \gamma_i\leq 0 \\
         &\Ab \xx_i \leq \bb,  \, \forall i=1,...,N,\\
         & \xx_i = \yy_j, \, \forall i \in \mathcal{E}_j, \, \mathcal{E}_j \in \mathcal{E} \label{consensus-primal}\\  
         & -{\tt M}_j\textcolor{black}{\delta_{jk}^q}\leq y_{jk} \leq {\tt M}_j\textcolor{black}{\delta_{jk}^q}, \, \forall j=1,\ldots,K, \, \forall k=1,\ldots,n.
\end{align}
\end{subequations}
Problem \eqref{dcco-lfc-sos1-bigm-primal} is a distributed convex NLP problem and its solution has the advantage of generating linearizations about points closer to the feasible region. Therefore, primal solutions and primal bounds are obtained by iteratively solving problem \eqref{dcco-lfc-sos1-bigm-primal}. 

In the case that the primal problem is a centralized NLP, a feasible point that satisfies all linear and nonlinear constraints is considered to be the primal solution candidate. However, when the primal problem has to be solved distributedly, as in problem \eqref{dcco-lfc-sos1-bigm-primal}, some numerical considerations have to be taken into account. Particularly, in this case, in addition to all linear and nonlinear constraints, the consensus constraints \eqref{consensus-primal} have to be satisfied which is more challenging to deal with since all computational nodes have to agree on a consensus solution. In case the primal solution does not satisfy the consensus constraints within an acceptable numerical tolerance, poor outer approximations are generated and added to the dual problem. Therefore a larger number of iterations are required by the distributed NLP solver, especially when a large computational network is considered. Another challenge in solving the primal problem distributedly is the communication burden between the nodes of the network and the LFCs.

\section{Distributed Hybrid Outer Approximation}\label{sec:algorithms}
In this section, we propose the Distributed Hybrid Outer Approximation (DiHOA) algorithm to solve problem \eqref{dcco-lfc-sos1-bigm-epig}. In essence, DiHOA is developed based on the LP/NLP-based BnB algorithm proposed in \cite{quesada1992lp}. The LP/NLP-based BnB algorithm is an implementation of the standard OA algorithm, where only a single BnB tree is built and outer approximations are added dynamically to the MIP master problem. Since only one BnB tree is constructed during the solution procedure, the LP/NLP-based BnB method is also called the single-tree OA algorithm. Similarly, the standard implementation of the OA algorithm where a BnB tree is constructed at each iteration of the algorithm is called multiple-tree OA. 

Before discussing the DiHOA algorithm, we briefly present Distributed Primal Outer Approximation (DiPOA) which is proposed in \cite{olama2021dipoa}. \vb{DiPOA} is the baseline algorithm developed to solve problem \eqref{dcco-lfc-sos1-bigm-epig} which, in essence, extends the standard OA algorithm so that the main computational parts related to the NLP sub-problems are handled in a distributed fashion.
From the numerical point of view, \vb{DiPOA} is a hierarchical algorithm which consists of three main computational layers, namely, \textit{primal}, \textit{cutting-plane manager}, and \textit{master} levels. 

The primal level deals with the nonlinear optimization part, particularly problem \eqref{dcco-lfc-sos1-bigm-primal}, consisting of two sub-levels responsible for a specific computational task. 
The first sub-level aims to simultaneously solve multiple local nonlinear optimization problems in a fully decentralized fashion. The local nodes then \textit{synchronously} communicate to the second sub-level, which is responsible for another phase of computations. The second sub-level is usually responsible for aggregating the solution of local NLP problems connected to each LFC by 
a series of unconstrained NLP problems. The solution of this level is then sent to the first sub-level until a consensus is reached.

The main purpose of the cutting-plane manager level is to generate and manage the OA linearizations, which are obtained around the generated feasible points provided by the primal level. This level is also responsible for adding the linearizations into the master's problem \eqref{dcco-lfc-sos1-bigm-master}. 

Finally, the master's level corresponds to solving problem \eqref{dcco-lfc-sos1-bigm-master} based on the cuts generated by the cutting-plane manager. The master MIP problem approximates the nonlinear functions of problem \eqref{dcco-lfc-sos1-bigm-epig}. As the number of cuts increases, this approximation improves until a good piece-wise outer approximator is achieved. The binary solution of the master level is then sent to the primal level, which, together with the nonlinear optimization problems, solves another set of NLP problems.

Generally speaking, when solving problem \eqref{dcco-lfc} with DiPOA, most of the total solution time is usually spent on solving the MIP master's problem. In such a scenario, the MIP problems are similar in consecutive iterations since they only differ by a few linear constraints. In particular, at iteration $q$ of DiPOA, a new feasible point is provided by the primal, around which a new linear approximation constraint is generated and added to the master's problem. In the next iteration, $q + 1$, the master's problem is reconstructed and solved from scratch. 

Based on what we proposed in \cite{olama2021dipoa} and according to \cite{quesada1992lp}, we develop the DiHOA algorithm to avoid constructing many similar MIP BnB trees. The main idea of DiHOA is to iteratively build a single branch-and-bound tree whereby the primal problem is solved distributedly, while dynamically updating the dual problem \eqref{dcco-lfc-sos1-bigm-master} without reconstructing the branch-and-bound tree.  However, the single-tree OA algorithm may lead to a large BnB and weaker approximations. To avoid that, DiHOA introduces several second-order outer approximations of the nonlinear functions to the root of the BnB tree through an event-triggered scheme, which leads to a tighter problem representation and a BnB tree with a fewer number of nodes. This procedure constructs the first MIP dual problem and initiates the BnB algorithm. During the BnB search, as soon as a new integer-feasible solution is found, the primal problem \eqref{dcco-lfc-sos1-bigm-primal} is distributedly solved to determine whether a \textit{lazy constraint} removing this integer-feasible point should be generated. By lazy constraints, we mean  cutting planes that are lazily added to the MIP model whenever an integer feasible solution is found. At this point, the RHADMM algorithm is applied, and the new primal information is distributed to the computational nodes of the network. Then the generated lazy constraint is added to the current node, and all open nodes of the BnB tree and the search continues. Therefore, it is not required to reconstruct the BnB tree as in multiple-tree algorithms and the same BnB tree can be used after adding new linearizations as lazy constraints. Finally, the algorithm is terminated until the MIP integer relaxation results in a feasible integer solution to the MINLP problem \eqref{dcco-lfc-sos1-bigm-epig}.

A detailed description of the DiHOA algorithm is summarized in Algorithm \ref{alg:dihoa}. It can be observed in Algorithm \ref{alg:dihoa} that DiHOA consists of three primary computational phases, namely, \vb{Initialization}, \vb{MultipleTreeSearch}, and \vb{SingleTreeSearch} steps. After the algorithm is successfully initialized, DiHOA starts a multiple-tree strategy with second-order outer approximations and continues the computations until either the solution is found or poor lower bound improvement is achieved. In the former case, the algorithm is terminated and the optimal solution is returned. In the latter case, however,  DiHOA accumulates all the outer approximations obtained until iteration $q$ in the root of the latest BnB tree and, then, it starts a single-tree search strategy whereby approximations are added dynamically. As the name of the algorithm suggests, DiHOA is a \textit{hybrid} algorithm that combines both single-tree and multiple-tree strategies by introducing an event-triggered scheme that determines the switching iteration,  $q_{switch}$, at which a single-tree search strategy is started. In the following, we describe each computational step in detail.

\subsection{Initialization Step}
The initialization step is started by solving the integer relaxation of problem \eqref{dcco-lfc-sos1-bigm-epig} to construct the dual problem \eqref{dcco-lfc-sos1-bigm-master} in the first iteration of the algorithm. The integer relaxation of problem \eqref{dcco-lfc-sos1-bigm-epig} is a consensus optimization problem that is solved using the RHADMM algorithm. In case the relaxation obtains a feasible solution with respect to problem \eqref{dcco-lfc-sos1-bigm-epig}, we terminate DiHOA with the optimal solution. Otherwise, we generate $N$ first-order outer approximations by using the local information available in each computational node of the network and construct the first MIP dual problem according to \eqref{dcco-lfc-sos1-bigm-master}.

\subsection{Multiple Tree Search Step}
The multiple-tree search step of the DiHOA algorithm is a crucial step that directly affects the DiHOA performance. In the cases that $q_{switch}$ is a large number, a pure multiple-tree algorithm with second-order outer approximations is obtained. Otherwise, the resulting algorithm becomes the single-tree OA. Therefore, $q_{switch}$ should be determined in such a way that maximizes the DiHOA performance. To do so, we introduce an event-triggered scheme that selects $q_{switch}$ based on the difference between two consecutive lower bounds (\textit{i.e.}, $lb^q - lb^{q-1}$). In particular, during the solution procedure, DiHOA checks if the generated lower bounds by the dual problem start to flatten out within a given tolerance $\epsilon > 0$ and triggers a switching event, $E_s$, if $lb^q - lb^{q-1} \leq \epsilon$. As soon as $E_s$ is triggered, DiHOA switches to the single-tree strategy by performing the BnB algorithm on the latest dual problem, which was obtained during the multiple-tree strategy.

\subsection{Single Tree Search Step}
The single-tree search step is activated for $q > q_{switch}$ after the $E_s$ event is triggered. In this step, DiHOA accumulates all the outer approximations obtained during the  \vb{MultipleTreeSearch} phase in the root of the latest BnB tree and then starts the single-tree BnB procedure. The BnB search is initialized by solving an integer relaxation of the dual problem \eqref{dcco-lfc-sos1-bigm-master}. In each node of the BnB tree, a QCLP relaxation is solved and the search is stopped once an integer solution is obtained in one of the nodes. The integer solution is then used to solve the primal problem \eqref{dcco-lfc-sos1-bigm-primal} with integer variables fixed. The primal solution provides a valid upper bound and new approximations can be generated. The new approximations are then added to all open nodes in the BnB tree and the QCLP relaxation is resolved for the node which resulted in the integer combination. The BB procedure continues from the existing search tree with the improved polyhedral outer approximation.

As in the standard BnB, nodes can be pruned off in case the optimum of the QCLP relaxation exceeds the upper bound. However, the search cannot be stopped once an integer solution is obtained at a node, which must continue until the QCLP relaxation results in a feasible integer solution for problem \eqref{dcco-lfc-sos1-bigm-epig} or until the node can be pruned off. Finally, since all variables in problem \eqref{dccp} are bounded, then the assumptions A1--A3 in \cite{Fletcher1996} and assumptions in \cite{Kronqvist2020} are valid and Slater's constraint qualification holds for problem \eqref{dcco-lfc-sos1-bigm-epig} when the binary variables are fixed. Hence, no NLP subproblem will be infeasible and the convergence of Algorithm \ref{alg:dihoa} is ensured. 

\begin{algorithm}
\SetAlgoLined
\caption{DiHOA Algorithm}\label{alg:dihoa}
\begin{enumerate}
    \item \vb{Initialization}
    \begin{enumerate}
        \item[1.1] obtain a relaxed solution $\Tilde{\xx}_i^0$, $\Tilde{\yy}_j^0$, and $\Tilde{\pmb{\delta}}_j^0$, $i \in \{1,\dots,N\}$,  $j \in \{1,\dots,K\}$ by solving an integer relaxation of the MINLP problem \eqref{dcco-lfc-sos1-bigm-epig}
        \item[1.2] generate $N$ outer approximations according to \eqref{foc} and construct the dual problem \eqref{dcco-lfc-sos1-bigm-master}
        \item[1.3] store the generated outer approximations.
        \item[1.4] set $q = 1$, $ub^0 = +\infty$, $lb^0 = -\infty$
    \end{enumerate}
    \item \vb{MultipleTreeSearch}
    \item[] \While{$lb^q - lb^{q-1} > \epsilon$}{
    \begin{enumerate}
        \item[2.1] solve problem \eqref{dcco-lfc-sos1-bigm-master} and obtain a dual solution $\pmb{\delta}_j^q$
        \item[2.2] $lb^q \xleftarrow{} \sum_{i=1}^{N}\gamma_i$
        \item[2.3] solve problem \eqref{dcco-lfc-sos1-bigm-primal} and obtain a primal solution  $\xx_i^q$, $\yy_j^q$
        \item[2.4] $ub^q \xleftarrow{} \min\left(\sum_{i=1}^{N}f_i(\xx_i^q), ub^{q-1}\right)$
        \item[2.5] If $ub^q - lb^q \leq \tau$, return $\xx_i^q$ as an optimal solution $\xx_i^*$ 
        \item[2.6] generate and store outer approximations according to \eqref{soc} and update dual problem \eqref{dcco-lfc-sos1-bigm-master} by adding new outer approximations
        \item[2.7] $q\xleftarrow{}q+1$
    \end{enumerate}
}
\item \vb{SingleTreeSearch}
\begin{enumerate}
  \item[3.1] start BnB search for the MIP problem obtained in the last step with all second-order outer approximations accumulated in the root
\end{enumerate}

\item[] \While{$ub^q - lb^{q} > \epsilon$}{
    \begin{enumerate}
        \item[3.2] if a new feasible integer solution, $\Bar{\pmb{\delta}}_j^q$, is found, then  $\pmb{\delta}_j^q \xleftarrow{} \Bar{\pmb{\delta}}_j^q$ and update $lb^q$ according to the current BnB tree
        \item[3.3] solve primal problem \eqref{dcco-lfc-sos1-bigm-primal} and update the primal solution  $\xx_i^q$, $\yy_j^q$
        \item[3.4] $ub^q \xleftarrow{} \min\left (\sum_{i=1}^{N}f_i(\xx_i^q), ub^{q-1}\right)$ $\xx_i^*$
        \item[3.5] generate outer approximations according to \eqref{foc} and add them to the current and open nodes of BnB
        \item[3.6] resolve the integer relaxation problem for the node which resulted in the integer combination
        \item[3.7] $q\xleftarrow{}q+1$ and continue exploring BnB tree
    \end{enumerate}
}
\end{enumerate}
\end{algorithm}

\section{SCOT: A Software Framework for Sparse Optimization}\label{sec:scot}
This section discusses the technical features of \vb{SCOT}, including architecture and design, main components, basic syntax, and usage. \vb{SCOT} is entirely written in \vb{C++17} programming language and at its core uses the \textit{Message Passing Interface (MPI)} \cite{gropp1999using} to perform distributed computation operations, communication and message-passing between nodes of the network. To handle local NLP optimization problems, \vb{SCOT} relies on various open source solvers, such as OSQP \cite{stellato2020osqp} and IPOPT \cite{wachter2006implementation}. \vb{SCOT} also implements a truncated Newton's method \cite{nocedal2006numerical} to solve unconstrained optimization problems. Moreover, \vb{SCOT} integrates commercial and open source MIP solvers such as \vb{Gurobi} \cite{gurobi}  and \vb{CBC} to solve MIP dual problems.  To solve the distributed NLP problem \eqref{dcco-lfc-sos1-bigm-primal}, \vb{SCOT} implements the RHADMM algorithm proposed in \cite{olama2019} using the main MPI operations. Finally, both DiPOA and DiHOA algorithms are implemented within \vb{SCOT} and are available through \vb{SCOT} Python API, \vb{ScotPy}, or \vb{SCOT} Command Line Interface (CLI).

\subsection{Architecture}
\begin{figure}
	\centering
	\includegraphics[width=0.6\textwidth]{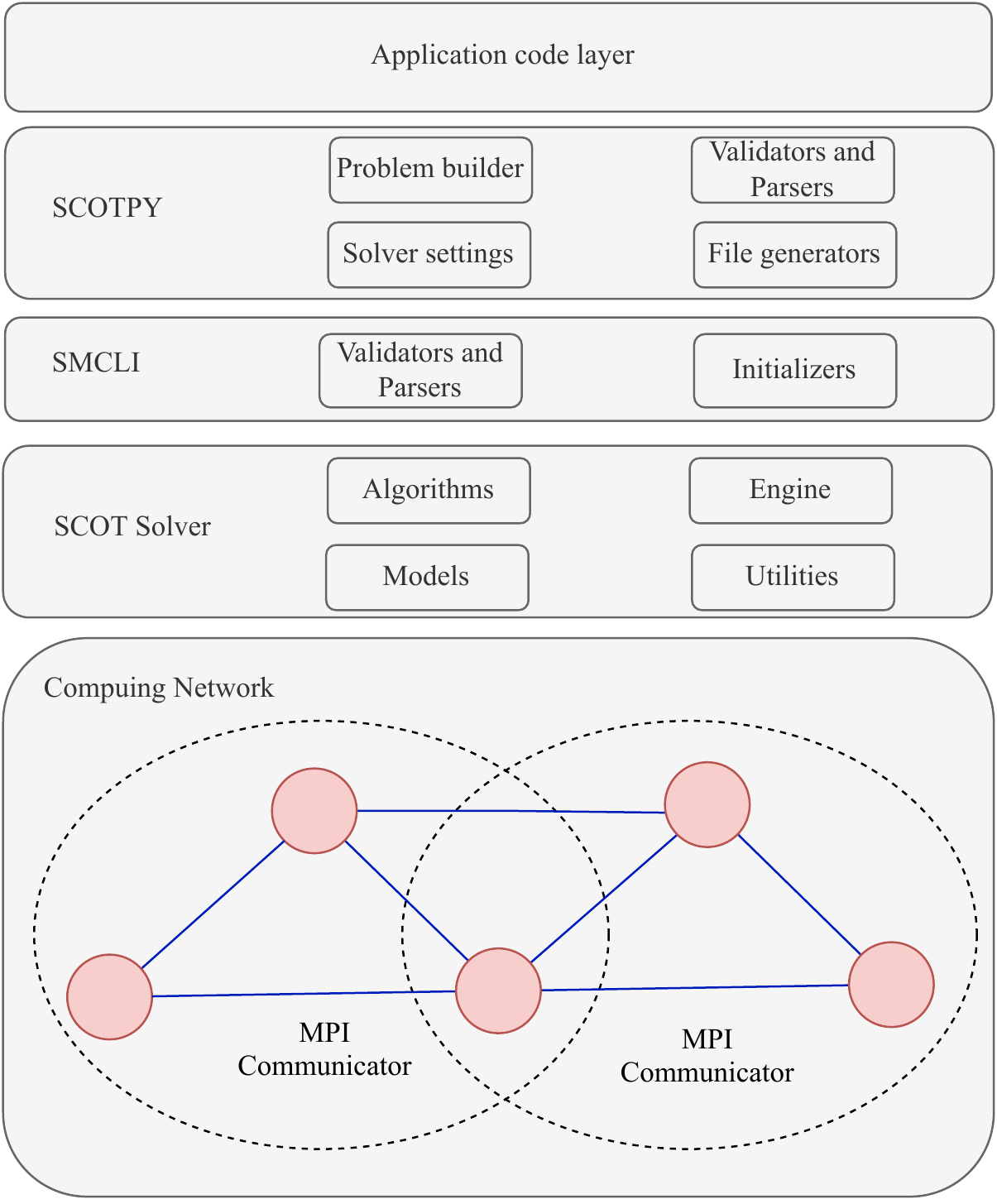}
	\caption{SCOT Architecture} 
	\label{sys-arch}
\end{figure}

A high-level overview of \vb{SCOT} architecture, its main layers, and components are shown in Figure \ref{sys-arch}. As observed in the figure, the primary layers of the framework are \vb{SCOTPY}, \vb{SMCLI}, \vb{SCOT Solver}, and \vb{Computing Network} each of which consists of different tools and components to build and solve the optimization problem. The layered architecture of \vb{SCOT} leads to a highly modular framework that can be easily extended by new features and algorithms. In the following, we describe each layer of  \vb{SCOT} and its components. 

\subsubsection{SCOTPY}
\vb{SCOTPY} is the main API written in \vb{Python 3.8} programming language, which provides various modules and classes to define the optimization problem and solver settings.  This layer consists of four components, namely, problem builder, validators and parsers, solver settings, and file generators. In principle, \vb{SCOTPY} receives the problem input and algorithm settings from the application code layer and, after parsing and validating the input data and settings, writes the optimization problem and settings in \vb{JSON} format with a specific naming convention. These steps are performed by validators and parsers, and file generator components. Moreover, according to the number of nodes given by the user, $N$ objective functions with different problem data are created and stored on the file system as different \vb{JSON} files. 

Therefore, the resulting computation of this layer is to express the optimization problem and solver settings in various \vb{JSON} files that can be read by the subsequent layers. Representing the optimization problem using files provides more flexibility since it decouples the optimization model from the optimization solver. Therefore, it is possible to write the optimization model in any language of choice and call the optimization solver from a different programming language or framework. Coupling the optimization model and the optimization solver through file formats is well-known in the optimization software industry and has been widely used for decades. We refer the interested reader to \cite{legat2021mathoptinterface} for more details. 

\subsubsection{SMCLI}
SCOT MPI Command Line Interface (SMCLI) is the main layer utilized to directly execute \vb{SCOT Solver} according to  different input and setting files. \vb{SMCLI} can be used directly without \vb{SCOTPY} interface, however the problem definition using \vb{SCOTPY} is more appropriate.  The validator and parsers component of \vb{SMCLI} is responsible to parse and validate the problem input and settings files, providing suitable data structures containing the optimization problem data for the \vb{SCOT Solver}. Additionally, initializing computational nodes and various software libraries used in \vb{SCOT Solver} are among the responsibilities of \vb{SMCLI} layer. 

\subsubsection{SCOT Solver}
At its core, the \vb{SCOT} framework consists of\vb{SCOT Solver} layer which is responsible to solve the optimization problem using a proper algorithm and settings. The main components of this layer are algorithms, engine, models, and utilities which are discussed in this section.

The \vb{Engine} component consists of various modules and classes to distributedly solve the primal problem  \eqref{dcco-lfc-sos1-bigm-primal} and deliver the primal solution to the MINLP algorithms.  Owing to its flexible and modular implementation, \vb{Engine} can be easily extended by introducing user-defined and custom-distributed convex optimization algorithms and solvers.    
In principle, \vb{Engine} requires the solution of local NLP problems for which multiple open-source and commercial solvers are available. At its core, \vb{Engine} consists of a sub-module, \vb{kernel} that provides a flexible interface to third-party solvers that can solve local optimization problems.
The solvers supported by the \vb{kernel} module are \vb{OSQP}, \vb{IPOPT}, and \vb{Gurobi}. As a final note, the \vb{Engine} component is called by all internal algorithms of \vb{SCOT} and handles most of \vb{MPI} communications and collective operations. More importantly, the quality of outer approximations depends on this component since it provides feasible points around which nonlinear functions are approximated. 

One of the most critical components of \vb{SCOT Solver} is \vb{algorithm} component that implements Algorithms \ref{alg:dihoa} and DiPOA.
The \vb{algorithms} component consists of two main modules, namely \vb{dipoa}, and \vb{dihoa}, which are responsible for implementing their corresponding algorithm.
Because all the implemented algorithms must manage and monitor outer approximations, the \vb{algorithm} component also includes a \vb{managers} module.
This module implements several classes to support the necessary data structures that generate and store both first- and second-order  outer approximations. 
Moreover, the \vb{managers} module implements the event-triggered schemes to improve the outer approximation's quality and switch from \vb{MultipleTreeSearch}  to \vb{SingleTreeSearch} strategy.
Among all classes, the \vb{managers} module consists of two important classes, namely \vb{CutStorage} and \vb{CutGenerator}. 
\vb{CutStorage} provides a simple way to validate and store the linear and quadratic outer approximations. The primary responsibility of the \vb{CutGenerator} class is to generate necessary outer approximations from the information received from the \vb{Engine} component. 
The \vb{algorithm} module also provides a flexible functionality to interface third-party MIP solvers such as \vb{Gurobi} and \vb{Cplex} for solving the MIP dual problem \eqref{dcco-lfc-sos1-bigm-master}.

The \vb{model} component's primary responsibility is to generate a concrete internal representation of the optimization problem to be used by algorithms. This component consists of various classes to present different types of nonlinear objective functions and linear and sparsity constraints. By accessing the \vb{model} component, the \vb{algorithm} will be able to access the optimization problem data whenever necessary during the computations.

Finally, the \vb{utilities} module implements classes and functions to provide commonly required functionalities, such as low-level parsing, measuring the CPU time, writing logs, constant parameters, exceptions, and file handling functions.

\subsubsection{Computing Network}
This layer is responsible for presenting and managing the computational network using a graph data structure and Message Passing Interface (MPI) library. The \vb{computing network} layer is in tight communication with the \vb{SCOT Solver} layer since all distributed algorithms use MPI for performing distributed computations and inter-process communications.

\subsection{Message Passing Interface (MPI)}
This section provides an overview of MPI and its important operations. \vb{MPI} is a message-passing library that supports parallel and distributed computations. For being independent of any programming language, MPI is arguably the most widely used platform for high-performance distributed computing nowadays. 
MPI is a software library for which various interfaces exist from various programming languages, including \vb{C/C++} and \vb{Python}.
MPI adopts the \textit{Single Program, Multiple Data} (SPMD) programming paradigm to provide an efficient way to perform distributed computing. 
   Using the SPMD paradigm, each node of the computing network runs the same program code, but it works with its own set of local variables and a separate subset of the data. 
To perform efficient distributed computing, MPI provides point-to-point and collective communications between the nodes of the computing network. 
   Point-to-point communications refer to sending and receiving data between two different nodes, whereas collective communication is primarily performed among a set of nodes. 
In the following, we review some of the standard collective operations.

\subsubsection{MPI Broadcast}
\begin{figure}
\centering
\subfloat[MPI broadcast.]{%
\resizebox*{7cm}{!}{\includegraphics{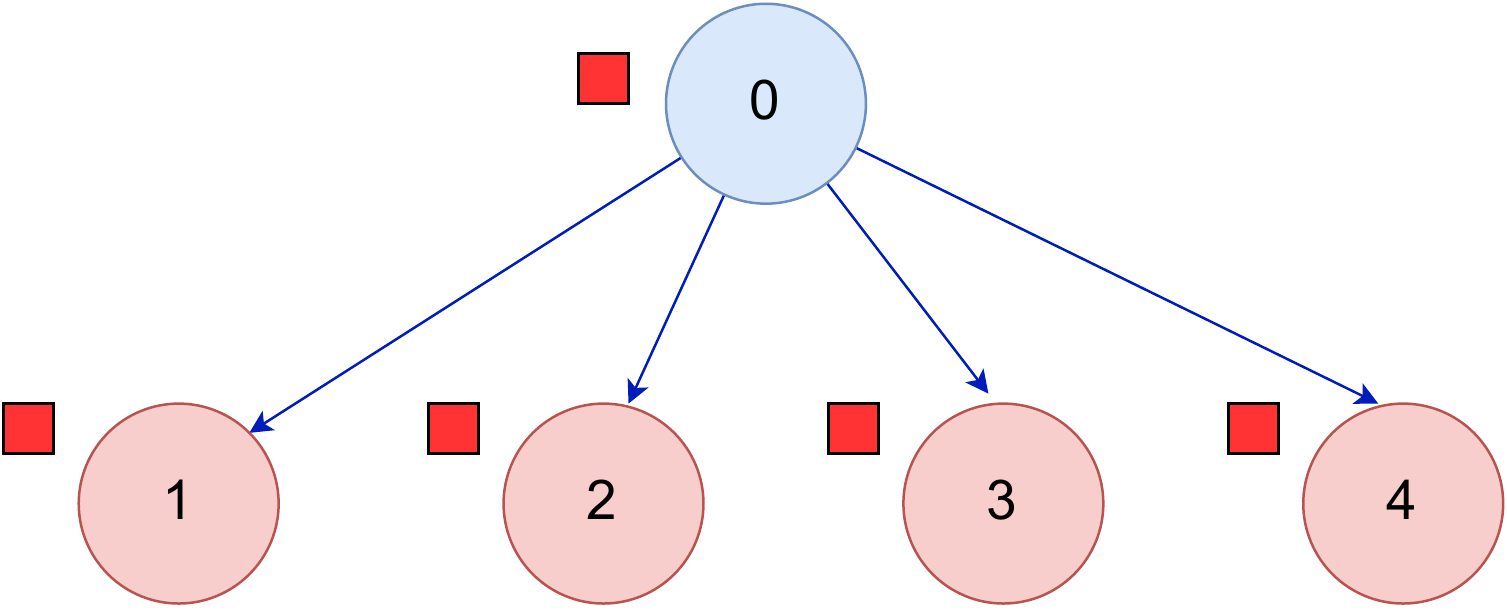}\label{bcast}}}\hspace{5pt}
\subfloat[MPI scatter.]{%
\resizebox*{7cm}{!}{\includegraphics{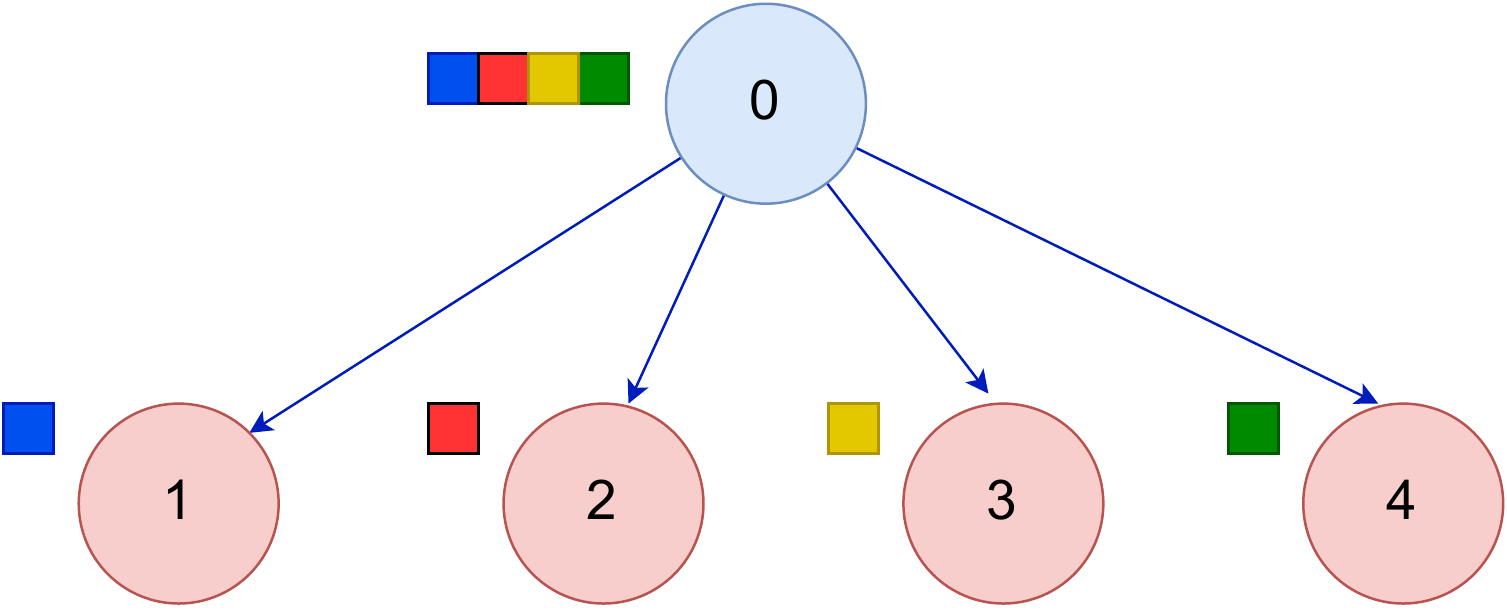}\label{scatter}}}	
\caption{MPI broadcast and scatter communication patterns.}  \label{mpi: bcast_scatter}
\end{figure}

A \textit{broadcast} is one of the basic collective communication strategies where one process sends the same data to all processes. Figure \ref{bcast} illustrates the  broadcast communication pattern.

\subsubsection{MPI Scatter}
A \textit{scatter} is a collective routine that involves a designated root process sending data to all other processes. 
  The main difference between MPI scatter and broadcast is that while the MPI broadcast sends the \textit{same} piece of data to all other processes, the MPI scatter sends \textit{chunks of an array} to different processes, as shown in Figure \ref{scatter}.
  
\subsubsection{MPI Gather and AllGather}
\begin{figure}
\centering
\subfloat[MPI Gather.]{%
\resizebox*{7cm}{!}{\includegraphics{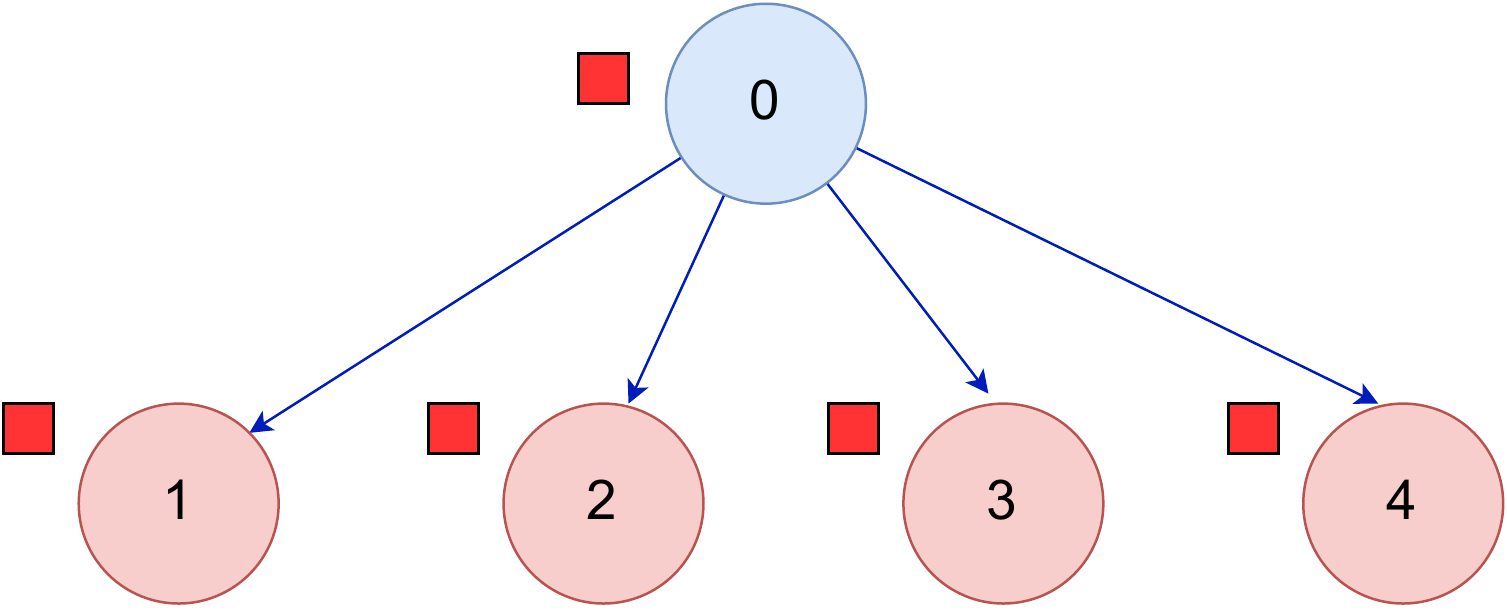}\label{gather}}}\hspace{5pt}
\subfloat[MPI Allgather.]{%
\resizebox*{7cm}{!}{\includegraphics{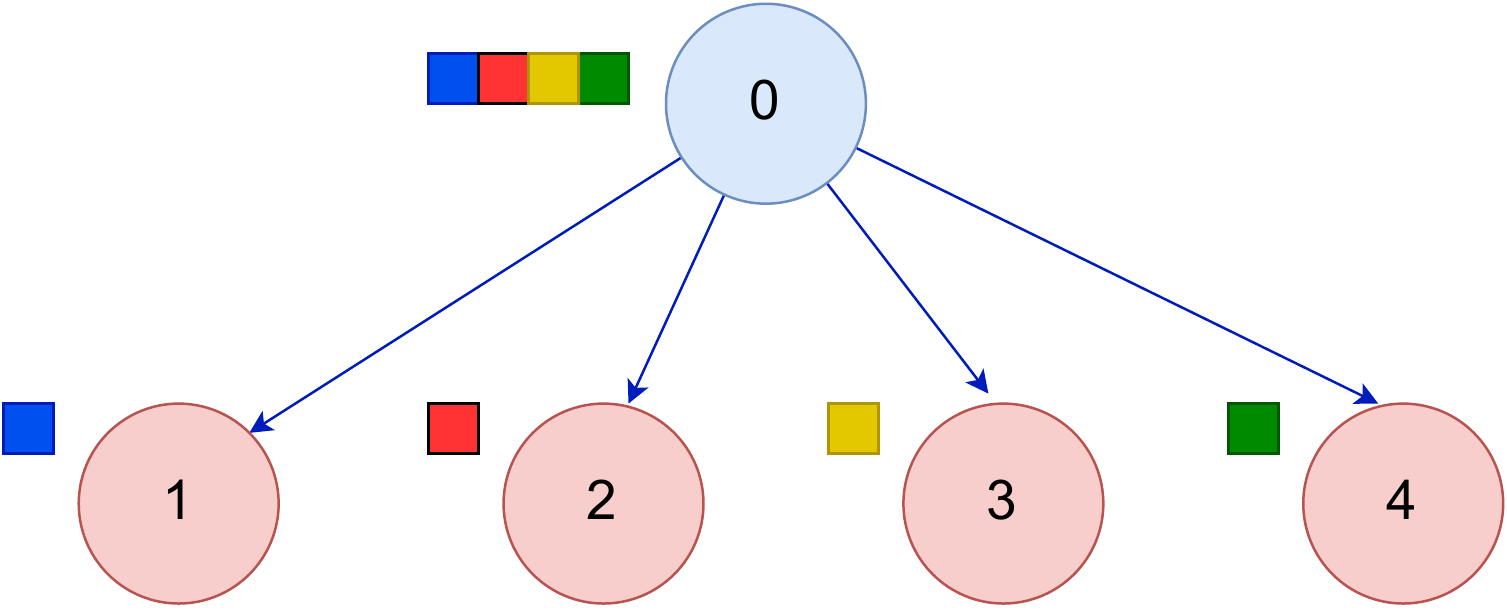}\label{allgather}}}
\caption{MPI Gather and Allgather communication patterns}  \label{mpi: gather_allgather}
\end{figure}

In principle, the MPI \textit{gather} is the inverse of the MPI scatter. Instead of distributing elements from one process to many processes, MPI gather takes elements from many processes and brings them together in one single process.
   The MPI \textit{Allgather} collects all of the elements and then distributes them to all the processes. 
In the most basic scenario, MPI allgather is an MPI gather followed by an MPI broadcast. For example, Figure \ref{mpi: gather_allgather} illustrates the communication pattern in MPI gather and allgather.

\subsubsection{MPI Reduce and Allreduce}
\begin{figure}
\centering
\subfloat[MPI Reduce.]{%
\resizebox*{7cm}{!}{\includegraphics{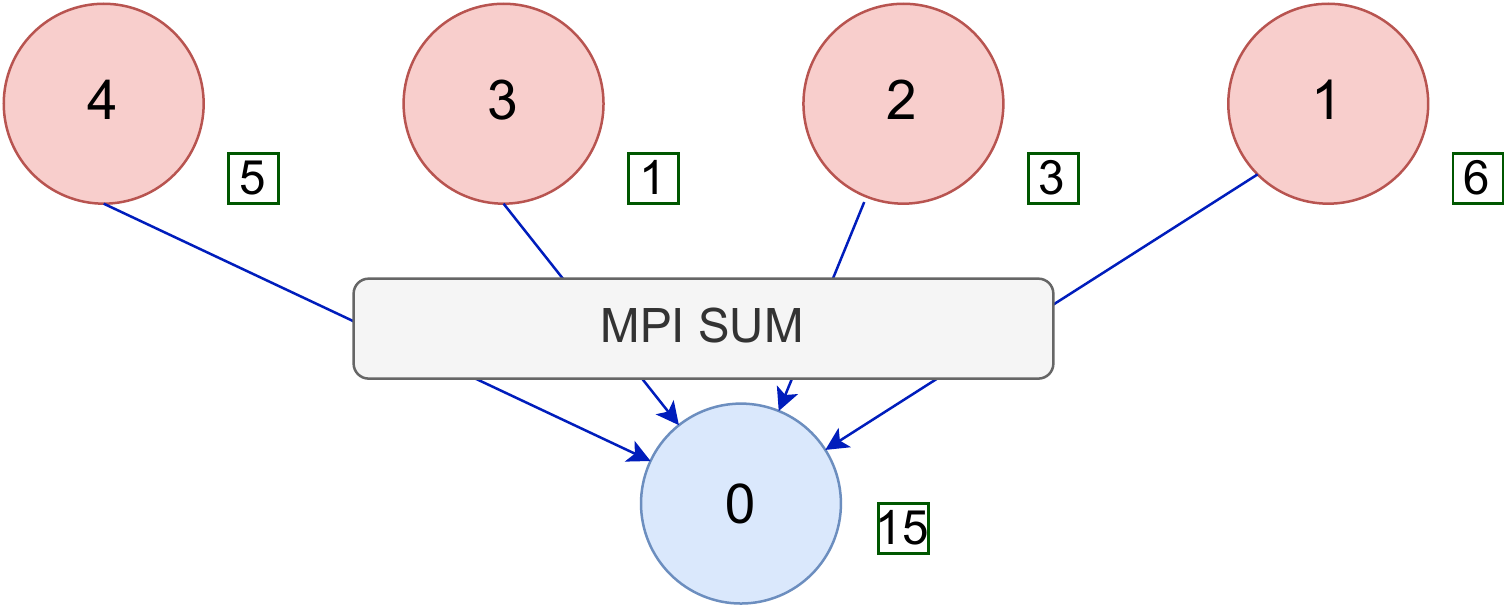}\label{reduce}}}\hspace{5pt}
\subfloat[MPI Allreduce.]{%
\resizebox*{7cm}{!}{\includegraphics{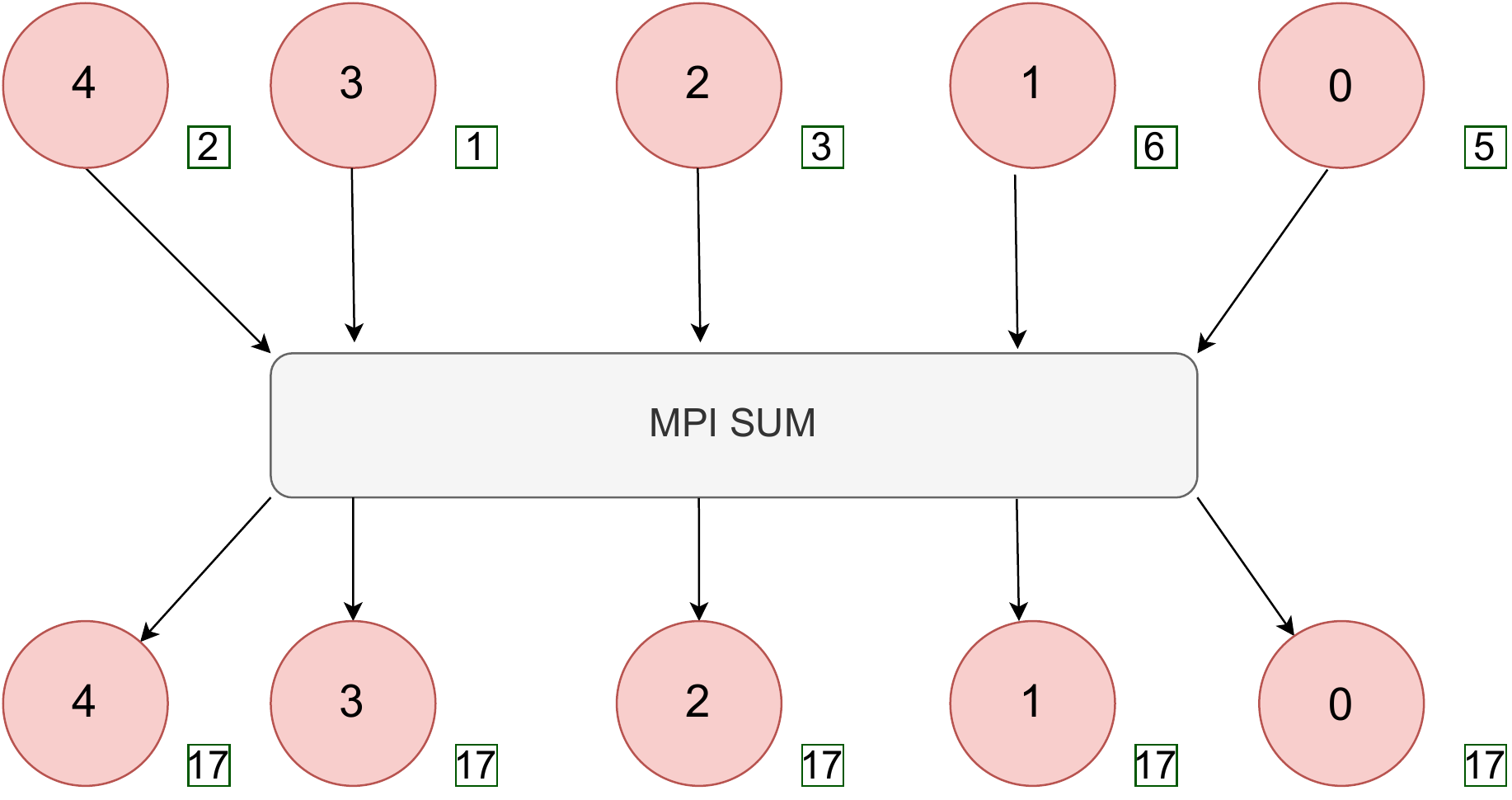}\label{allreduce}}}
\caption{MPI Reduce and Allreduce Communication Patterns.}  \label{mpi: reduce_allreduce}
\end{figure}

MPI \textit{reduce} involves reducing a set of elements into a small set of elements via a function. 
   The MPI reduce takes an array of input elements from each process and returns an array of output elements to the root process. 
 In a complementary style of MPI allgather to MPI gather, MPI allreduce will reduce the values and distribute the results to all processes. MPI reduce and allreduce communication patterns are depicted in Figure \ref{mpi: reduce_allreduce}. 

\subsection{Basic Syntax and Usage}
In this section, we present an illustrative example to show how \vb{SCOT} Python API, \vb{SCOTPY}, is used to solve a distributed sparse logistic regression problem with random data. To do so, we first import the required classes from \vb{SCOTPY} as the following code snippet shows,
\begin{lstlisting}[language=Python, caption= Import statement of SCOT Python API]

from scotpy import (AlgorithmType,
                    ProblemType,
                    ScotModel,
                    ScotPy,
                    ScotSettings
                    )
\end{lstlisting}

Here \vb{ScotPy} is the main class that executes \vb{SCOT} for a given problem  and settings defined by \vb{ScotModel} and \vb{ScotSettings}, respectively. The \vb{AlgorithmType} and \vb{ProblemType} classes determine what problem class is solved and which algorithm will be used. In order to create the optimization problem and \vb{SCOT} settings, the following code snippet can be used,
\begin{lstlisting}[language=Python, caption= Problem definition and settings]
# Create a classification dataset with 1000 rows and 20 columns
dataset, res = make_classification(n_samples = 1000, n_features = 20) 
scp = ScotModel(problem_name = "logistic_regression",
                rank = 0,
                kappa = 5,
                ptype = ProblemType.CLASSIFICATION)

# Set problem data with normalization
scp.set_data(dataset, res, normalized_data = True)

# Create corresponding JSON files that represent the optimization problem.
scp.create()

scot_settings = ScotSettings(
                relative_gap = 1e-5,
                time_limit = 100,
                verbose = True,
                algorithm = AlgorithmType.DIHOA
            )
\end{lstlisting}
where \vb{make\char`_classification} function, imported from Python \vb{scikit-learn} library, is used to generate a random classification dataset. The \vb{ScotModel} object is then created by a given problem name, MPI rank, number of nonzeros, and problem type. The solver settings can be defined by creating an object from \vb{ScotSettings} class. We note that by executing MPI, the above code snippets are simultaneously executed by each node of the network. Hence, each node can use its own problem data. Finally, we solve the optimization problem by using the following code,

\begin{lstlisting}[language=Python, caption= SCOT Execution]

solver = ScotPy(problem, scot_settings)
status_code = solver.run()

\end{lstlisting}

where \vb{ScotPy} class is responsible for creating a solver object for a given problem and settings. By executing the \vb{run} method of \vb{ScotPy}, MPI execution with $N$ nodes is started.


\section{Applications and Numerical Examples}\label{sec:applications}
In this section, we evaluate the \vb{SCOT} performance by comparing \vb{SCOT} to various state-of-the-art MINLP solvers equipped with single-tree, multiple-tree, and nonlinear BnB algorithms.  We considered \vb{SHOT} and \vb{BONMIN} as decomposition-based and \vb{KNITRO} as nonlinear BnB solvers. However, one should keep in mind that these are general-purpose solvers, and unlike \vb{SCOT} they are not tailored for this specific problem structure. From the application point of view, we focus on sparse logistic and linear regression problems with data distributed over a network of computational nodes. The benchmark results are obtained by generating performance profiles according to various problem instances over different computational nodes generated by \vb{SCOTPY}. 

\subsection{Implementation Details and Set-up}
All the experiments were performed on a Linux machine with an Intel Core i5 2.50 GHz processor, with four physical cores and 16 GB of RAM. The \vb{SCOT} solver is entirely implemented in \vb{C++} programming language and relies on \vb{MPI} to carry out distributed computation and inter-process communication. The source code of the solver is available on \url{https://github.com/Alirezalm/scot}. To perform linear algebra operations required by the distributed NLP solver, \vb{SCOT} uses \vb{Eigen} 3.4 library. Moreover, the MIP solver employed in \vb{SCOT} is \vb{GUROBI} 9.5.2 with an academic license. As for the comparison with other MINLP solvers, \vb{GAMS} 36.2 was selected as the optimization platform. It should also be noted that to achieve a meaningful comparison, \vb{GUROBI} is selected as the primary MIP solver for all MINLP solvers considered in the benchmarks.

\subsection{Experimental Set-up}
This section presents the problem classes that compose the benchmarks and experiments. We consider two well-known machine-learning problems: classification and regression. For both problem classes, the dataset is assumed to be distributed over a computational network, and the sparsity of the solution matters. Sparse classification and regression problems are among the central problems in statistics and machine learning, whose solutions usually lead to more interpretable models. In this paper, we choose logistic regression and linear regression as candidate models for classification and regression. 

Sparse classification and regression problems are tightly connected to SCO problems as it is often desired to identify a critical subset of features contributing to the response. Furthermore, the sparse solution usually leads to more interpretable models and improves prediction accuracy by eliminating unnecessary features. Accordingly, we introduce Distributed Sparse Logistic Regression (DSLogR) and Distributed Sparse Linear Regression (DSLinR) problems. The DSLogR problem is defined as,
\begin{subequations}\label{dslogr}
\begin{align}
      \min_{\thh_i, \yy_j, \delta_j} & \sum_{i = 1}^{N} \sum_{\ell = 0}^p \log\left[1 + \mathrm{e}^{-(\thh_i^T\X_{i,\ell})\Gam_{i,\ell}}\right] + \frac{\lambda}{2}\norm{\thh_i}_2^2\\
      & \thh_i = \yy_j, \, \forall i \in \mathcal{E}_j, \, \mathcal{E}_j \in \mathcal{E}, \, \forall j=1,...,K\\
       & -{\tt M}_j\delta_{jk} \leq y_{jk} \leq {\tt M}_j \delta_{jk}, \, \forall j=1,\dots,K, \, \forall k=1,\dots,n  \\
     &\textcolor{black}{\sum_{k=1}^{n} \delta_{jk} \leq \kappa, \, \forall j=1,\dots,K} \\
         & y_{jk}( 1 - \delta_{jk}^q) = 0, \, \forall j=1,\ldots,K, \, \forall k=1,\ldots,n\\
     & \delta_{jk} \in \{0, 1\}, \, \forall j=1,\ldots,K, \, \forall k=1,\dots,n 
\end{align}
\end{subequations}
where $\X_i \in \mathcal{R}^{ p \times n}$ and $\Gam_i\in \mathcal{R}^{ p}$ are the dataset and response vector of the $i$-th node, $\X_{i,\ell}$ is the column representation of $\ell$-th row of $\X_i$, $\Gamma_{i,\ell}$ it the $\ell$-th element of $\Gam_i$, and $\lambda >0$ is the regularization parameter. Clearly, problem \eqref{dslogr} is a subclass of the SCO problem \eqref{dcco-lfc-sos1-bigm-epig} which is the primary reformulation that \vb{SCOT} uses. Similarly, the DSLinR problem is defined as,
\begin{subequations}\label{dslinr}
\begin{align}
      \min_{\thh_i, \yy_j, \delta_j} & \sum_{i = 1}^{N} \norm{\mathbf{X}_i\thh_i - \mathbf{b}_i}_2^2 +  \frac{\lambda}{2}\norm{\thh_i}_2^2\\
      & \thh_i = \yy_j, \, \forall i \in \mathcal{E}_j, \, \mathcal{E}_j \in \mathcal{E}, \, \forall j=1,...,K\\
              & -{\tt M}_j\delta_{jk} \leq y_{jk} \leq {\tt M}_j \delta_{jk}, \, \forall j=1,\dots,K, \, \forall k=1,\dots,n  \\
         &\textcolor{black}{\sum_{k=1}^{n} \delta_{jk} \leq \kappa, \, \forall j=1,\dots,K} \\
          & y_{jk}( 1 - \delta_{jk}^q) = 0, \, \forall j=1,\ldots,K, \, \forall k=1,\ldots,n\\
         & \delta_{jk} \in \{0, 1\}, \, \forall j=1,\ldots,K, \, \forall k=1,\dots,n 
\end{align}
\end{subequations}
We generate $N$ random local datasets for both DSLogR and DSLinR problems with zero mean and unit $\ell_2$ norm for each column. 
We perform the numerical benchmarks based on different solver settings. Default settings were adopted for \vb{Knitro} as an NLP BnB solver. The chosen settings for \vb{SCOT}, \vb{SHOT}, and \vb{Bonmin} are reported in Table \ref{table:settings}.

\begin{table}
    \tbl{MINLP Solver Settings}
    {\begin{tabular}{lllll} \toprule
     & \multicolumn{2}{l}{Settings} \\ \cmidrule{2-5}
      solver & algorithm & name & mip solver & nlp solver  \\ \midrule
     {SCOT}  & {dipoa} & {scot-mt} &{gurobi} & {rhadmm} \\
     {SCOT}  & {dihoa} & {scot-st} &{gurobi} & {rhadmm}\\
     {SHOT}  & {ESH multiple-tree} & {shot-mt} &{gurobi} & {ipopt}\\
     {SHOT}  & {ESH single-tree} & {shot-st} &{gurobi} & {ipopt}\\
     {Bonmin}  & {B-OA} & {bonmin-mt} &{gurobi} & {ipopt}\\
     {Bonmin}  & {B-QG} & {bonmin-st} &{gurobi} & {ipopt}\\ \bottomrule
    \end{tabular}}
    \label{table:settings}
\end{table}

\subsection{Benchmark Results}

\begin{figure}
    \centering
    \subfloat[$90\%$ sparsity.]{%
    \resizebox*{7cm}{!}{\includegraphics{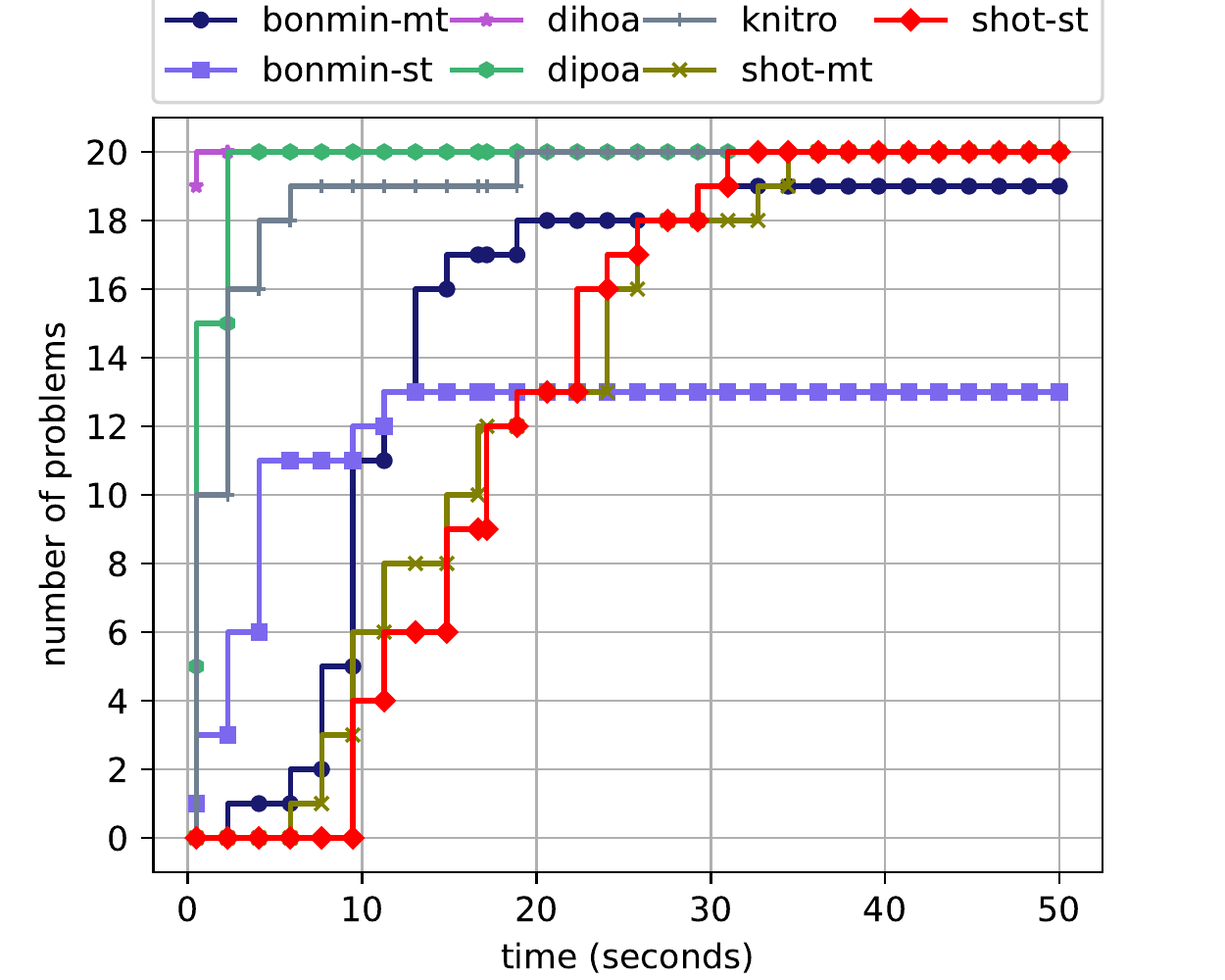}\label{fig:bench11}}}\hspace{5pt}
    \subfloat[$80\%$ sparsity.]{%
    \resizebox*{7cm}{!}{\includegraphics{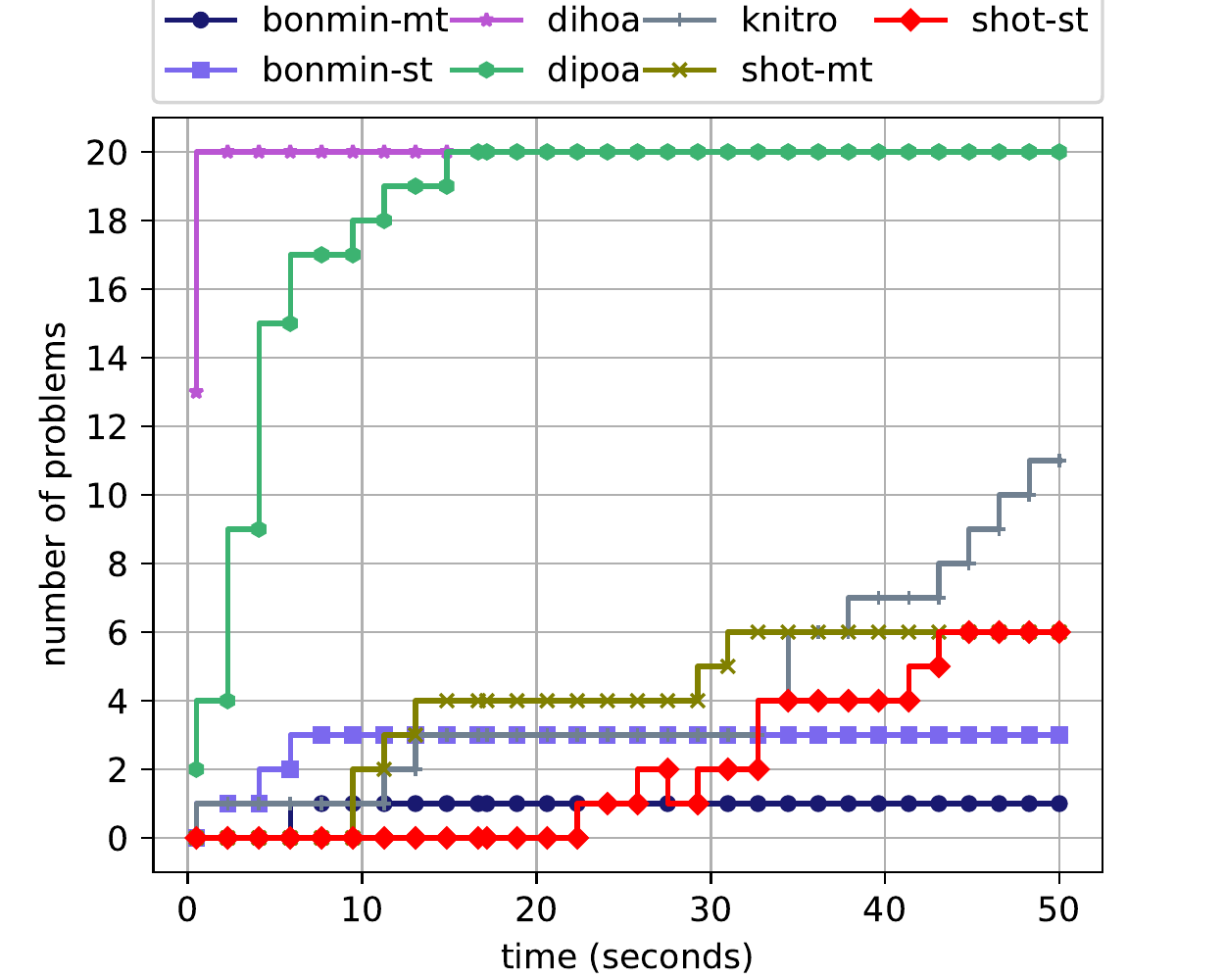}\label{fig:bench12}}}	
    \caption{Benchmark results for scenario $1$.}  \label{fig:bench1}
\end{figure}

\begin{figure}
    \centering
    \subfloat[$90\%$ sparsity.]{%
    \resizebox*{7cm}{!}{\includegraphics{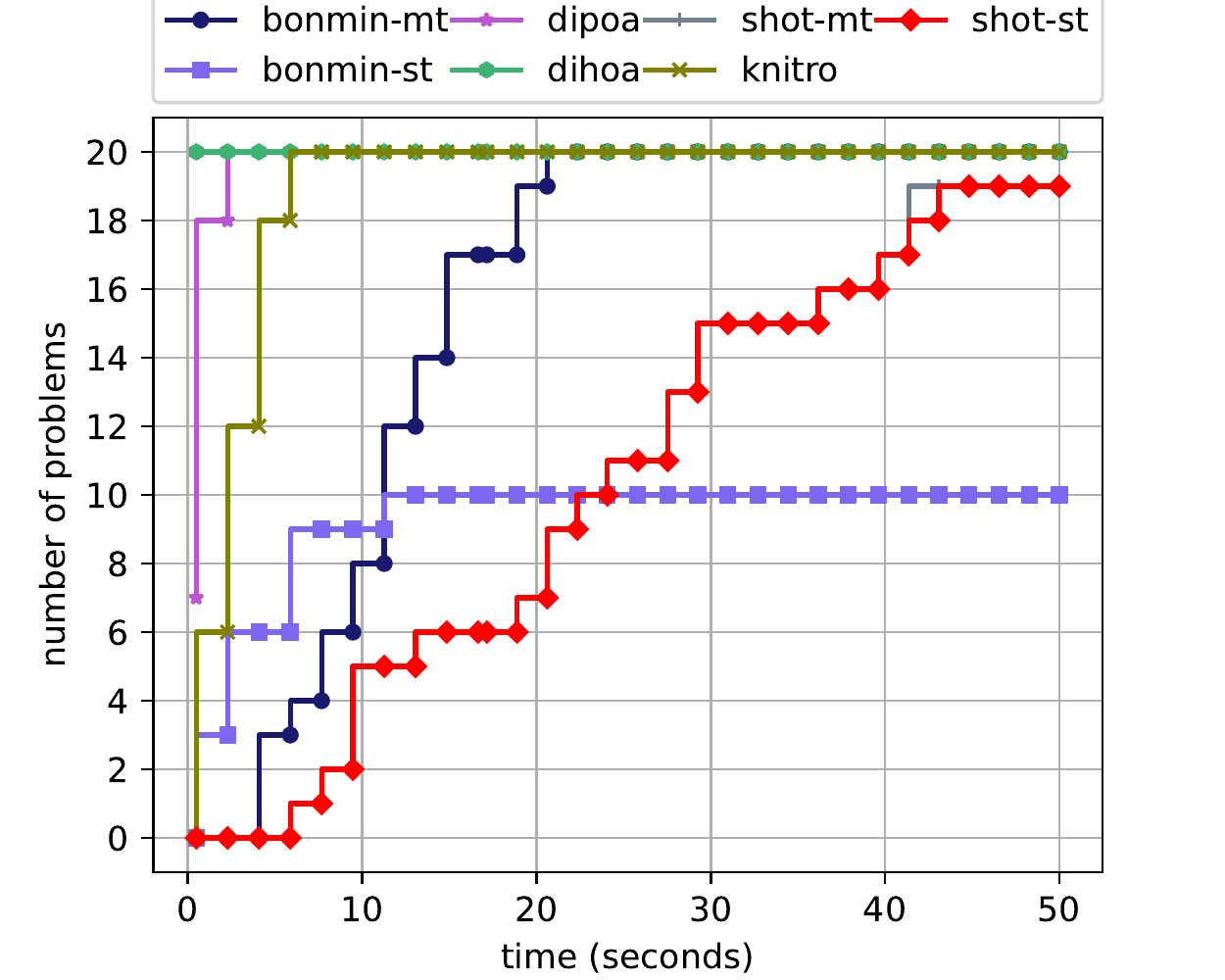}\label{fig:bench21}}}\hspace{5pt}
    \subfloat[$80\%$ sparsity.]{%
    \resizebox*{7cm}{!}{\includegraphics{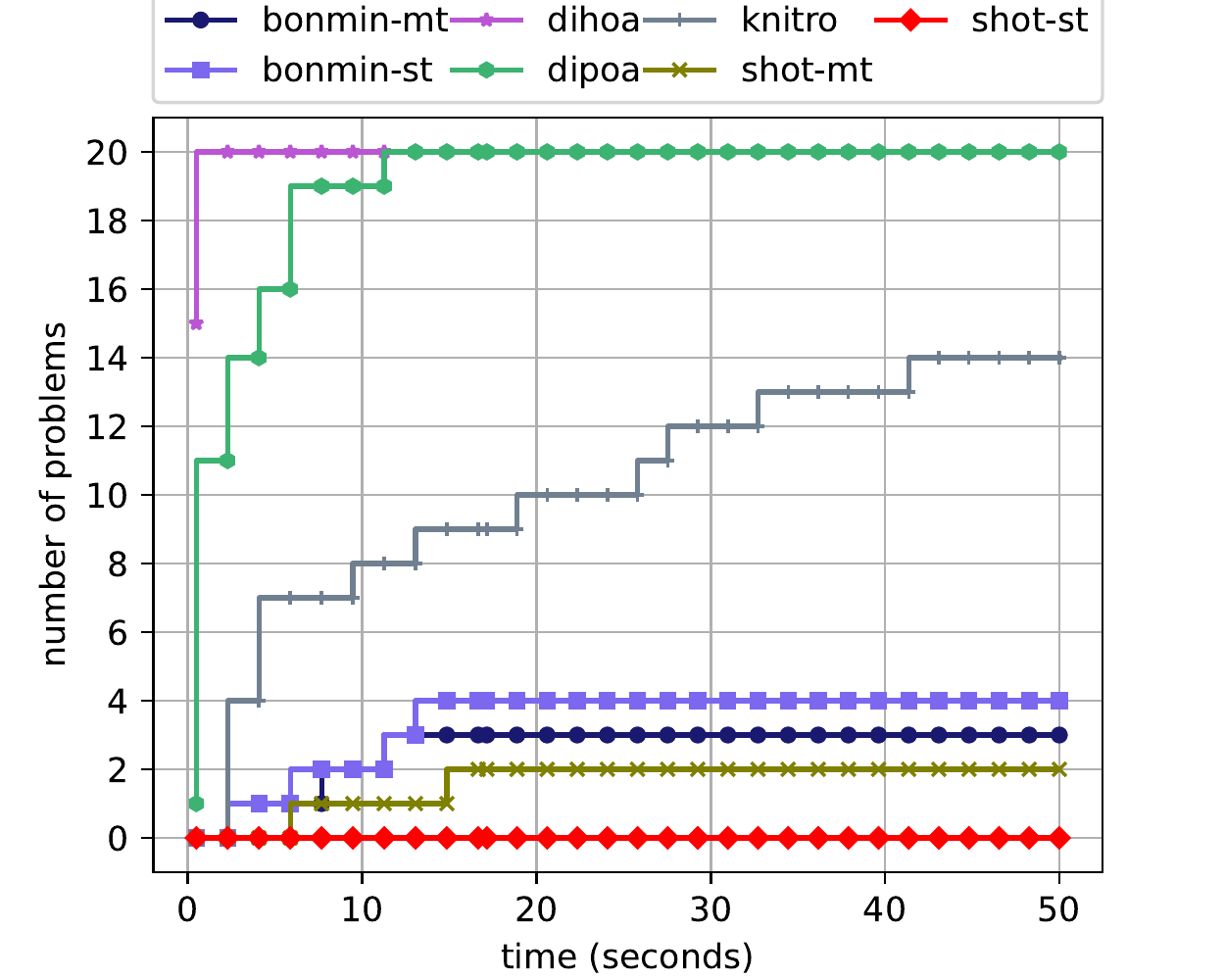}\label{fig:bench22}}}	
    \caption{Benchmark results for scenario $2$.}  \label{fig:bench2}
\end{figure}

\begin{figure}
    \centering
    \subfloat[$90\%$ sparsity.]{%
    \resizebox*{7cm}{!}{\includegraphics{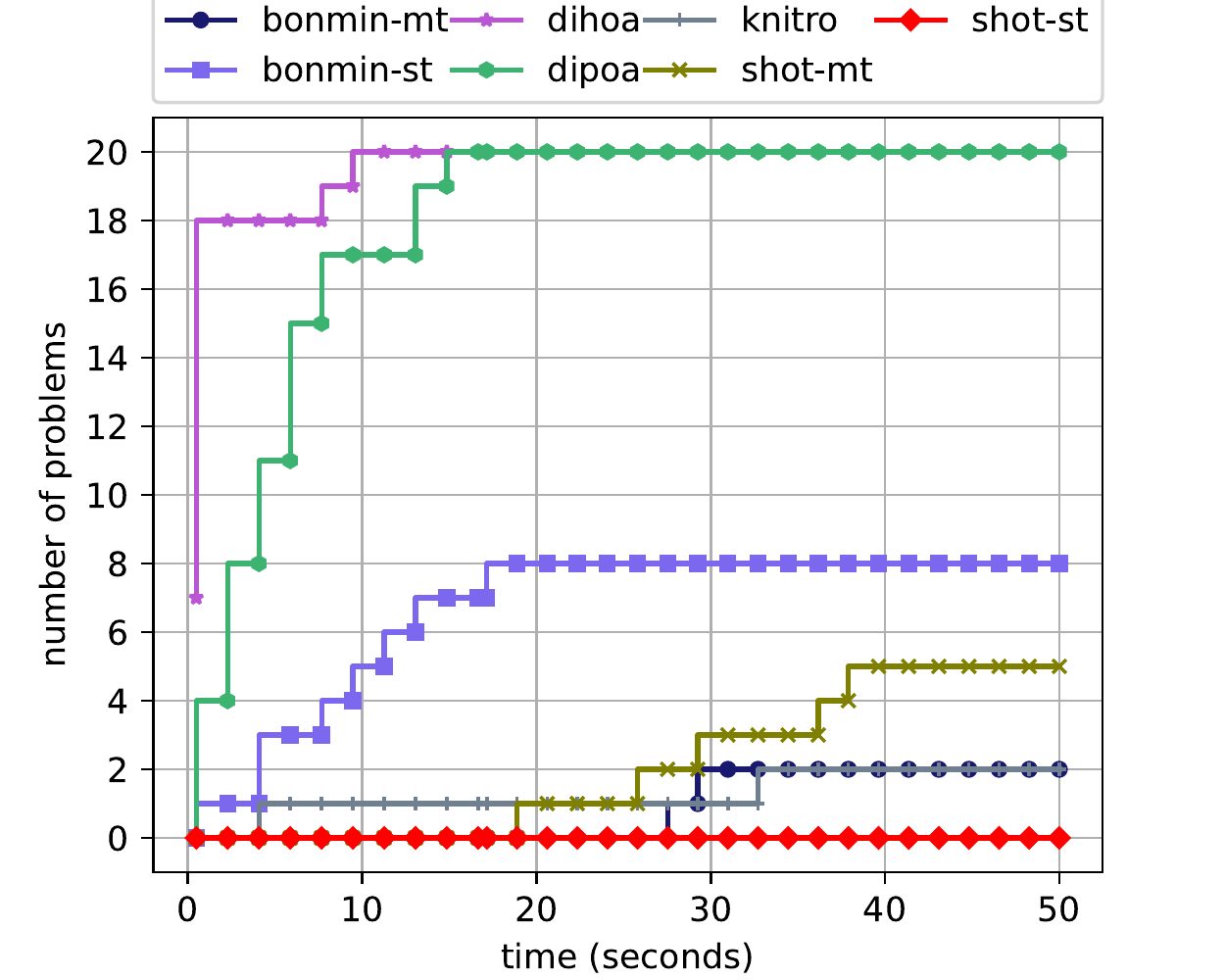}\label{fig:bench31}}}\hspace{5pt}
    \subfloat[$80\%$ sparsity.]{%
    \resizebox*{7cm}{!}{\includegraphics{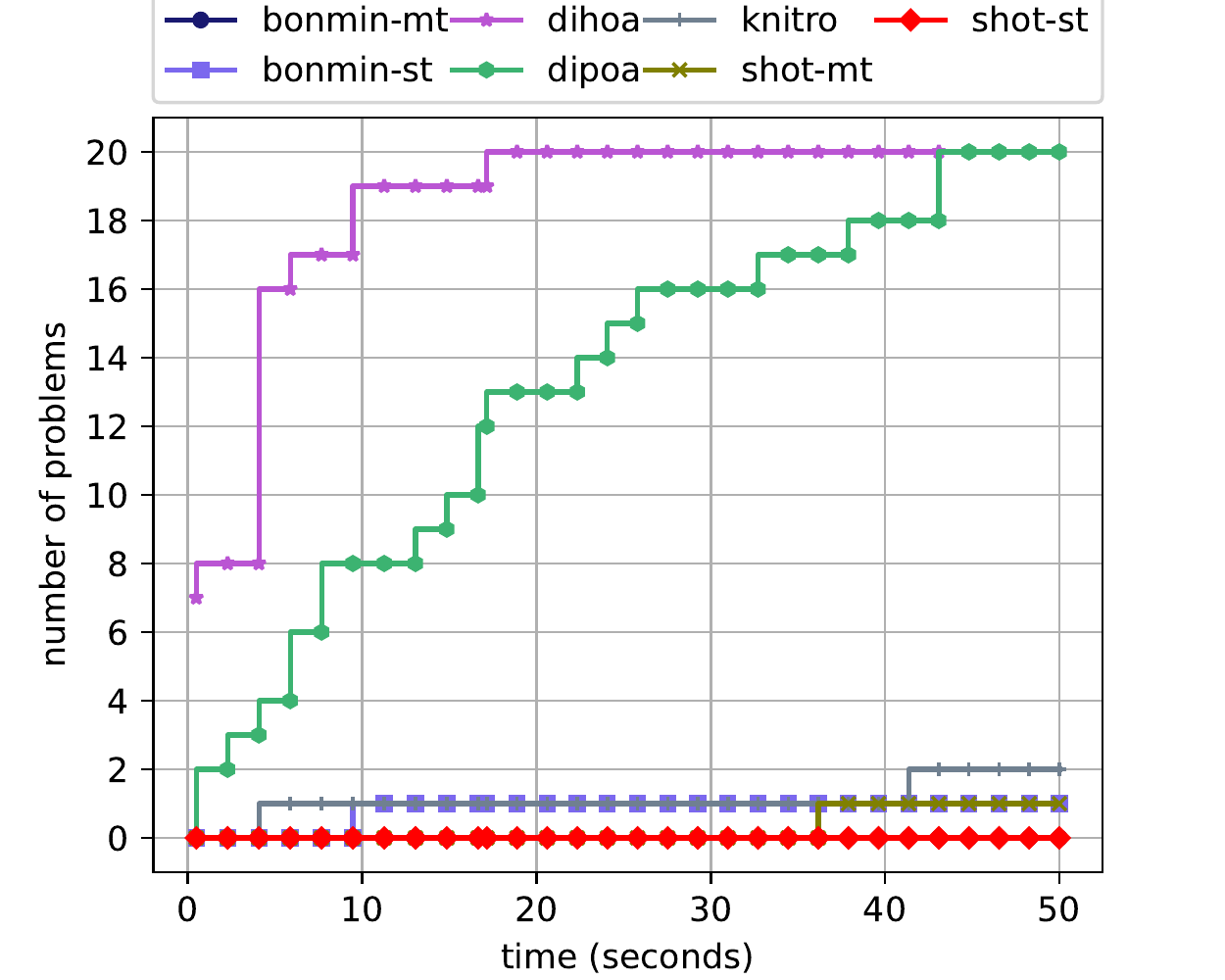}\label{fig:bench32}}}	
    \caption{Benchmark results for scenario $3$.}  \label{fig:bench3}
\end{figure}

\begin{figure}
    \centering
    \subfloat[$90\%$ sparsity.]{%
    \resizebox*{7cm}{!}{\includegraphics{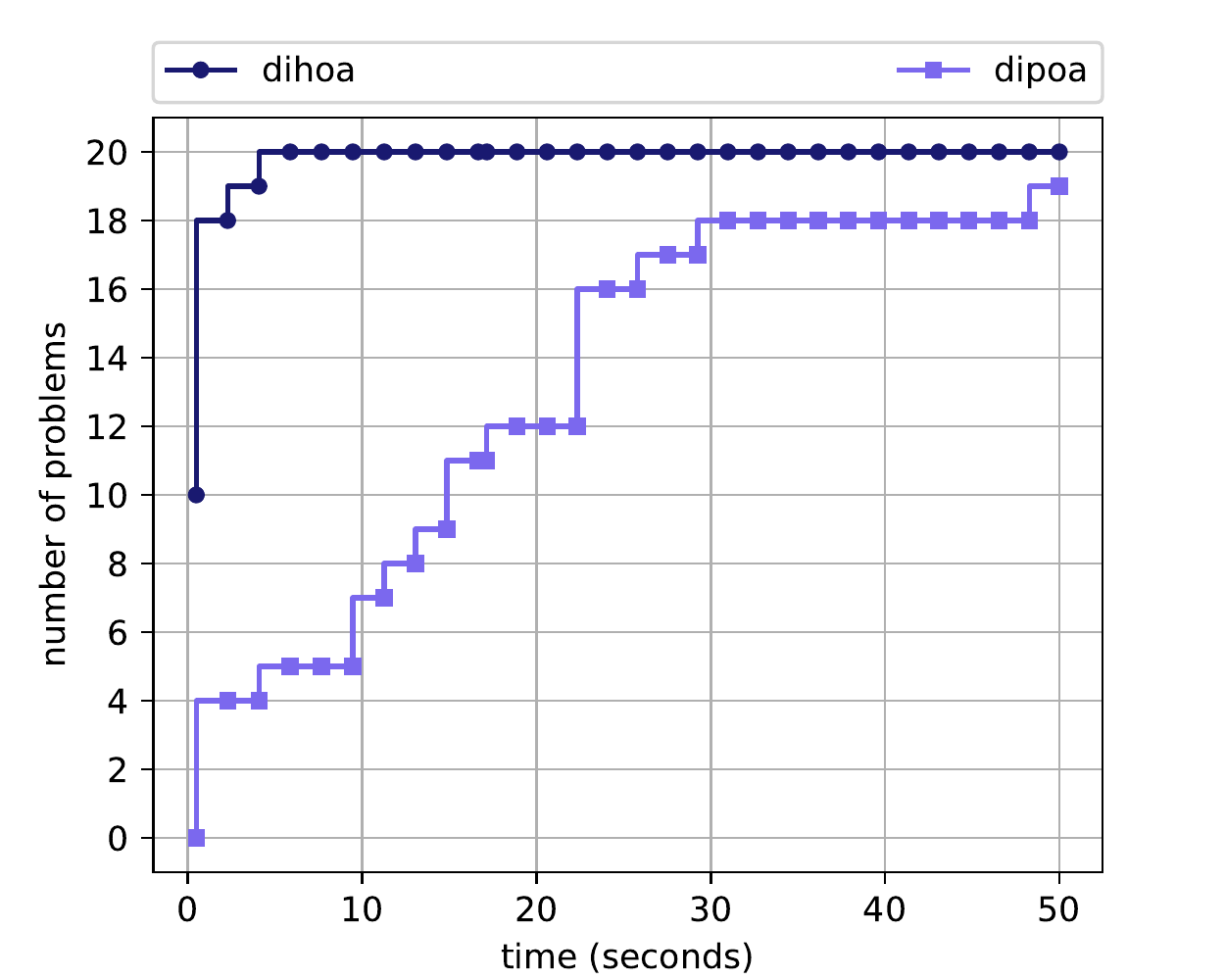}\label{fig:bench41}}}\hspace{5pt}
    \subfloat[$80\%$ sparsity.]{%
    \resizebox*{7cm}{!}{\includegraphics{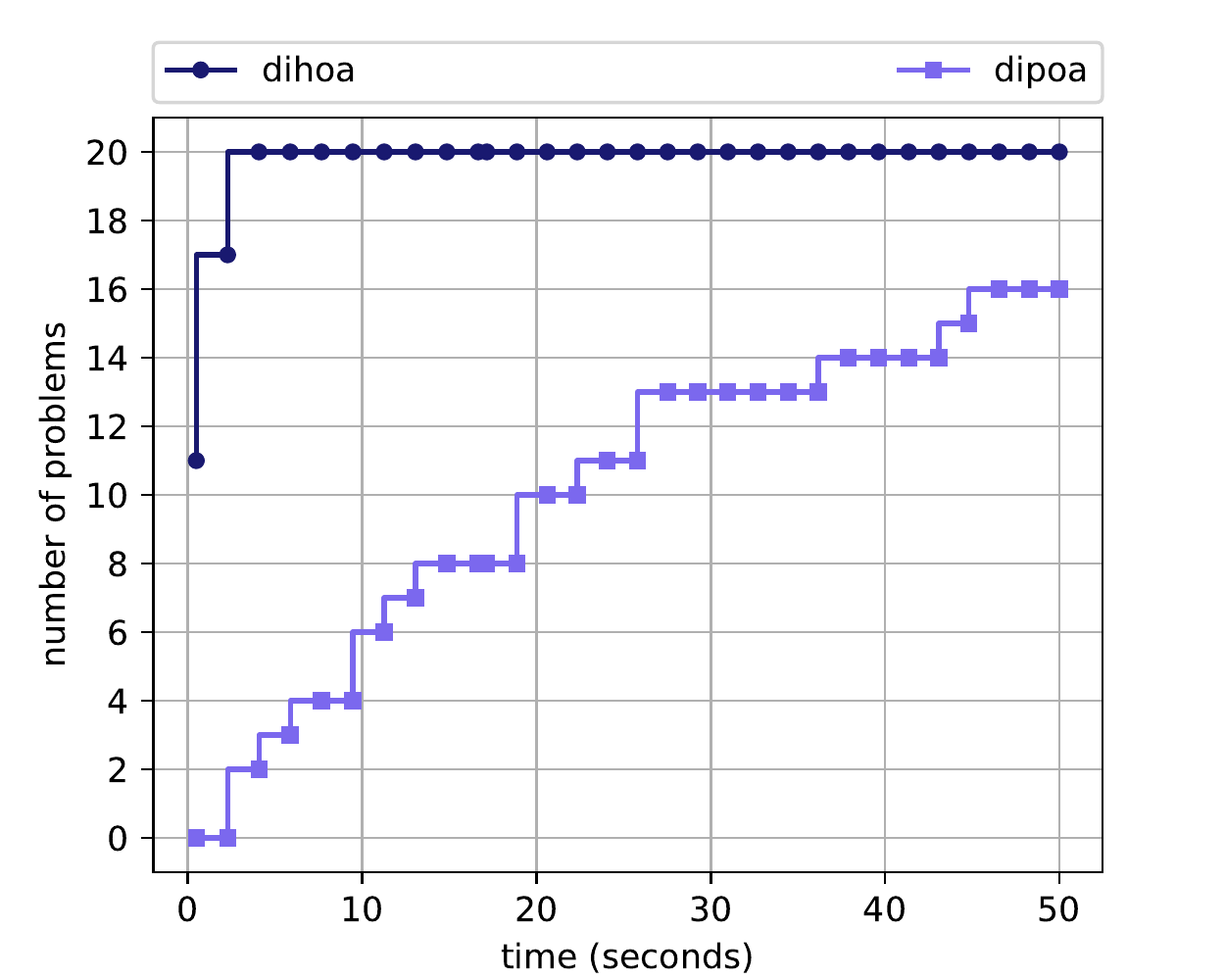}\label{fig:bench42}}}	
    \caption{Benchmark results for scenario $4$.}  \label{fig:bench4}
\end{figure}

\begin{figure}
    \centering
    \subfloat[$90\%$ sparsity.]{%
    \resizebox*{7cm}{!}{\includegraphics{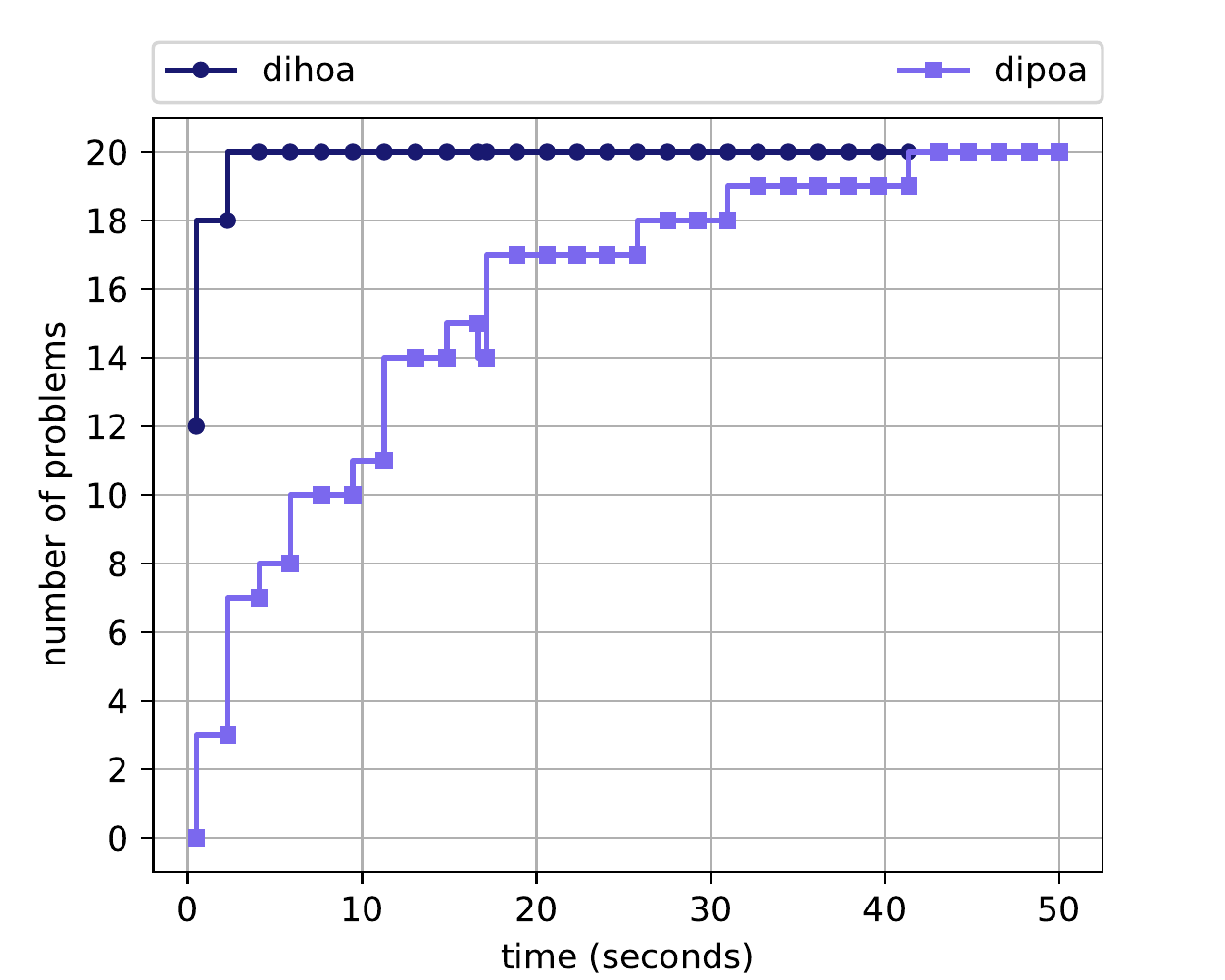}\label{fig:bench51}}}\hspace{5pt}
    \subfloat[$80\%$ sparsity.]{%
    \resizebox*{7cm}{!}{\includegraphics{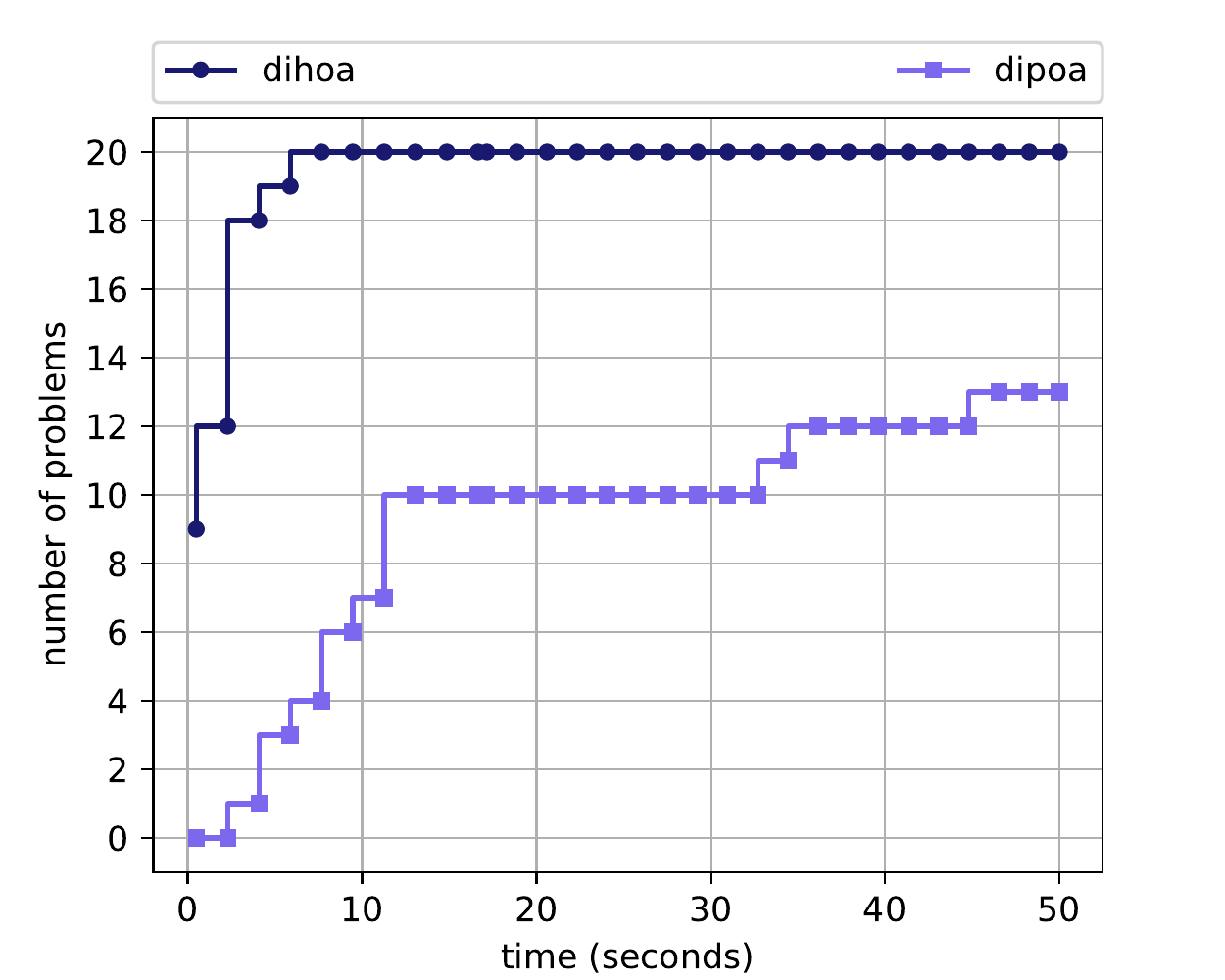}\label{fig:bench52}}}	
    \caption{Benchmark results for scenario $5$.}  \label{fig:bench5}
\end{figure}

\begin{figure}
    \centering
    \subfloat[$90\%$ sparsity.]{%
    \resizebox*{7cm}{!}{\includegraphics{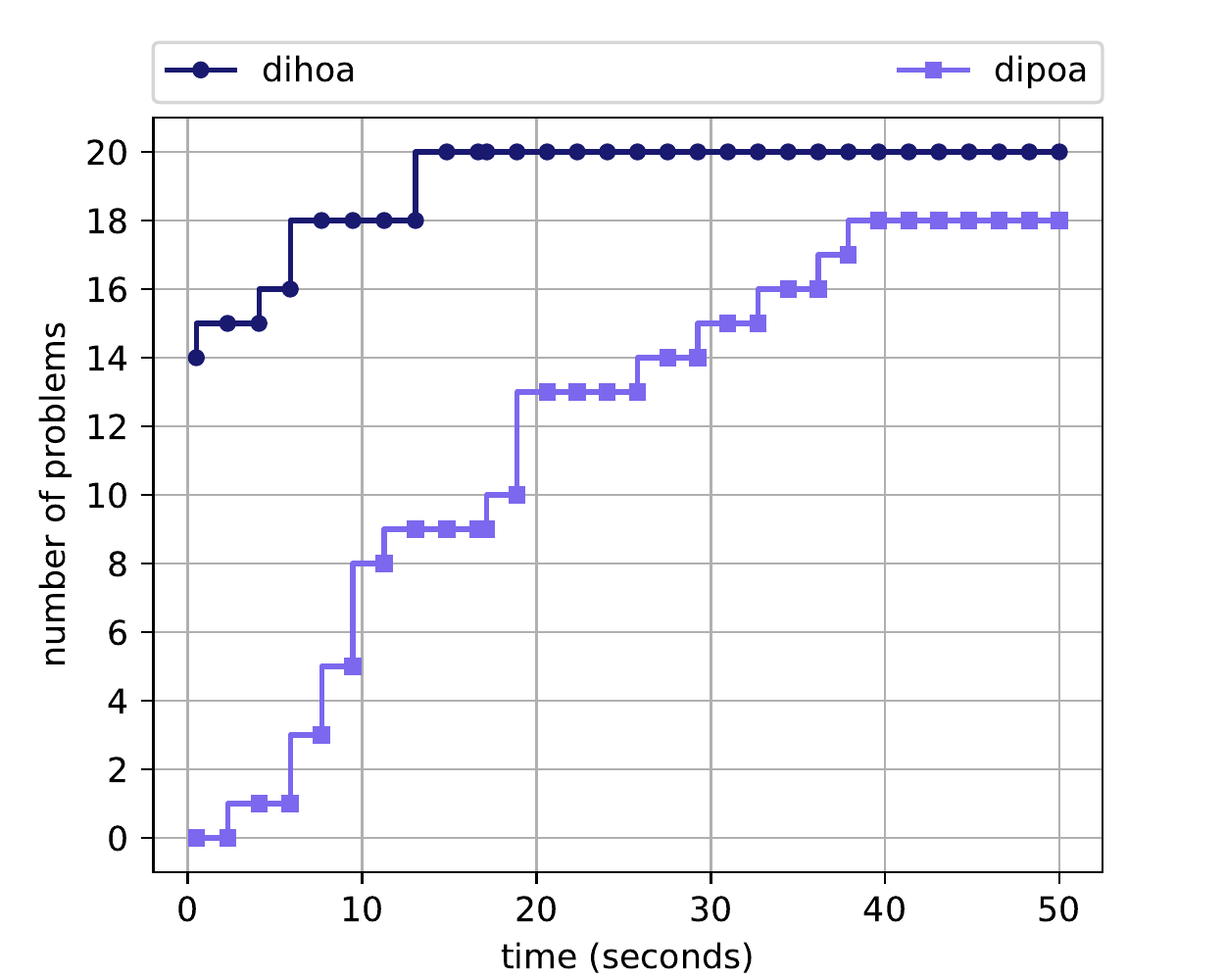}\label{fig:bench61}}}\hspace{5pt}
    \subfloat[$80\%$ sparsity.]{%
    \resizebox*{7cm}{!}{\includegraphics{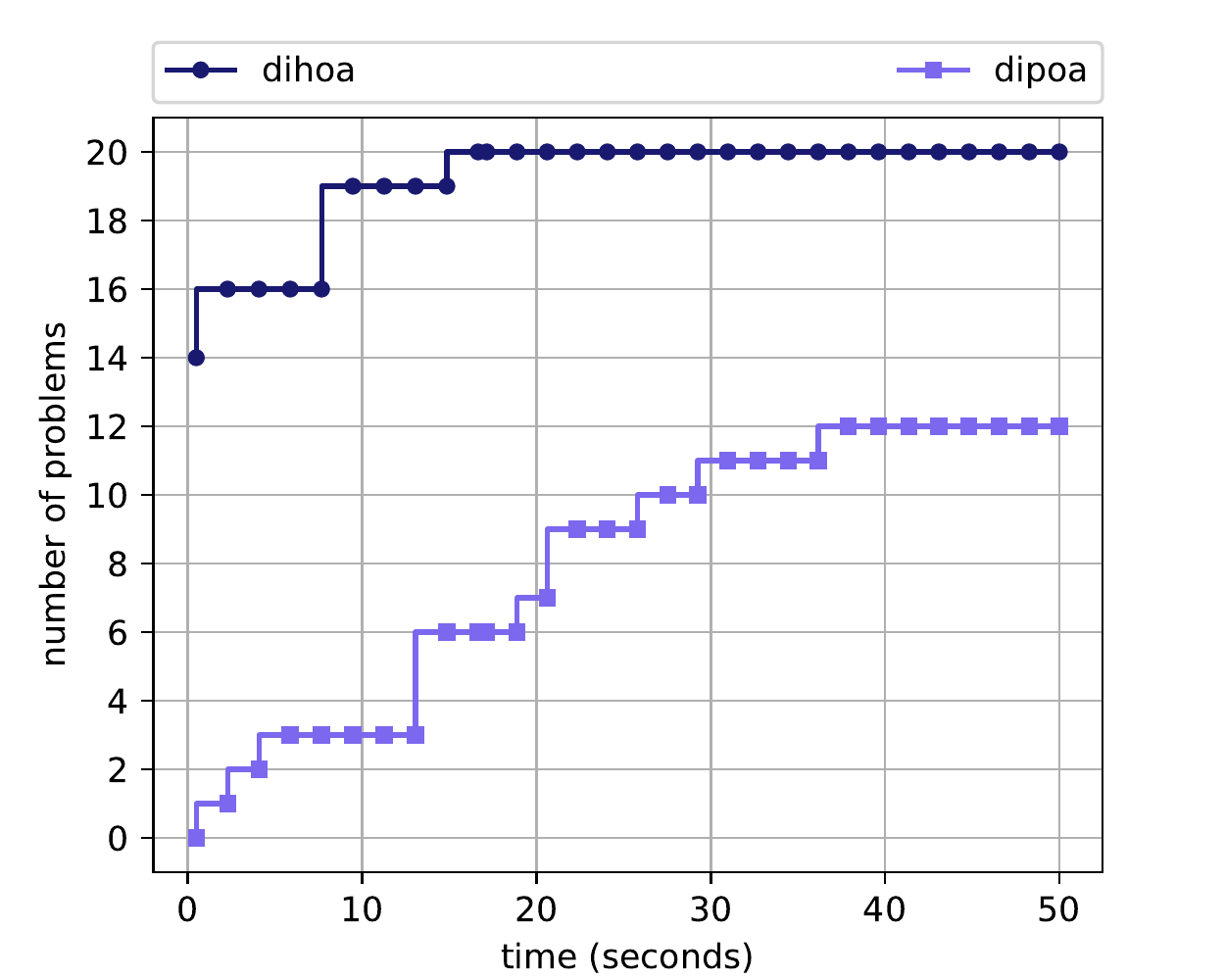}\label{fig:bench62}}}	
    \caption{Benchmark results for scenario $6$.}  \label{fig:bench6}
\end{figure}

\begin{figure}
    \centering
    \subfloat[$90\%$ sparsity.]{%
    \resizebox*{7cm}{!}{\includegraphics{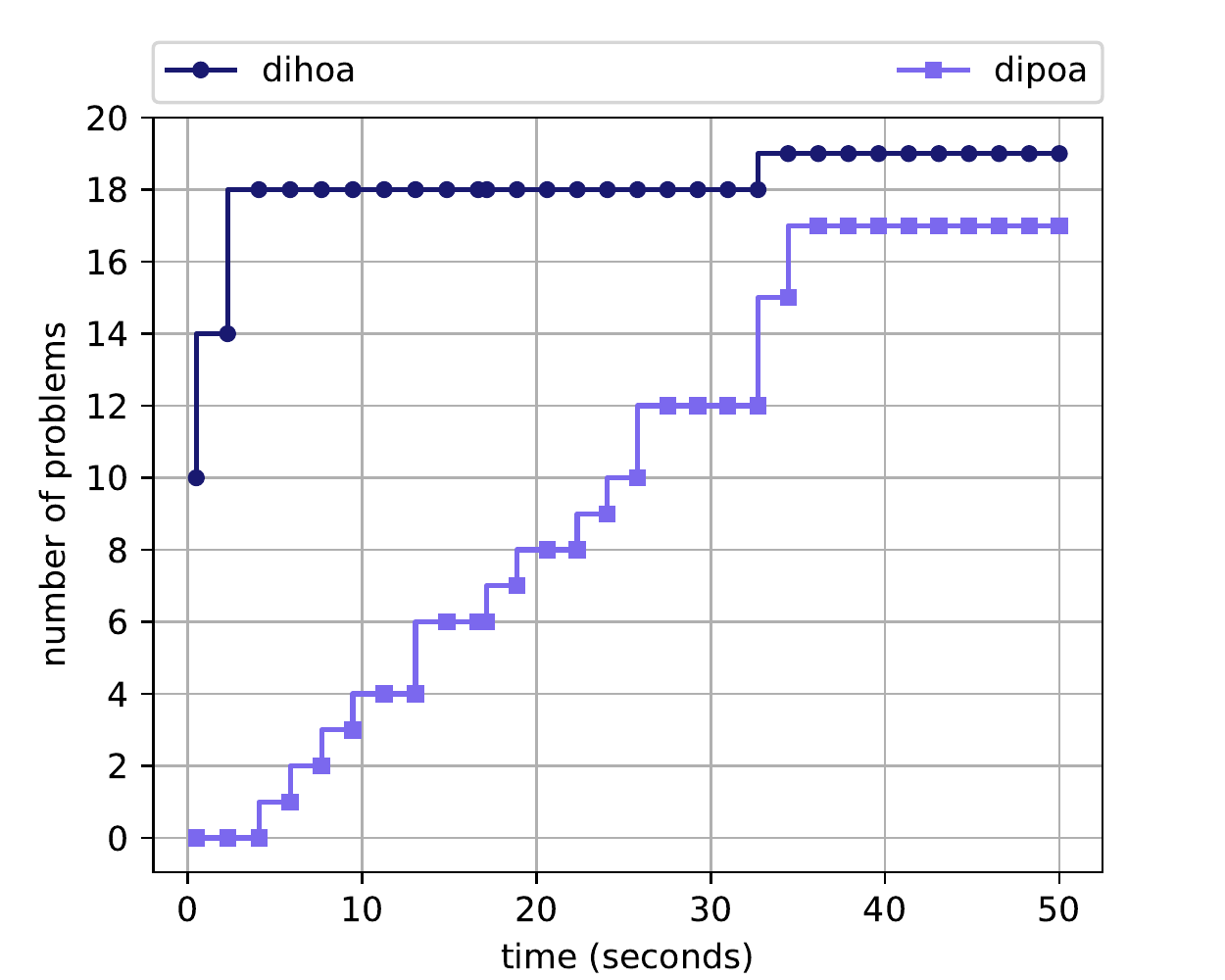}\label{fig:bench71}}}\hspace{5pt}
    \subfloat[$80\%$ sparsity.]{%
    \resizebox*{7cm}{!}{\includegraphics{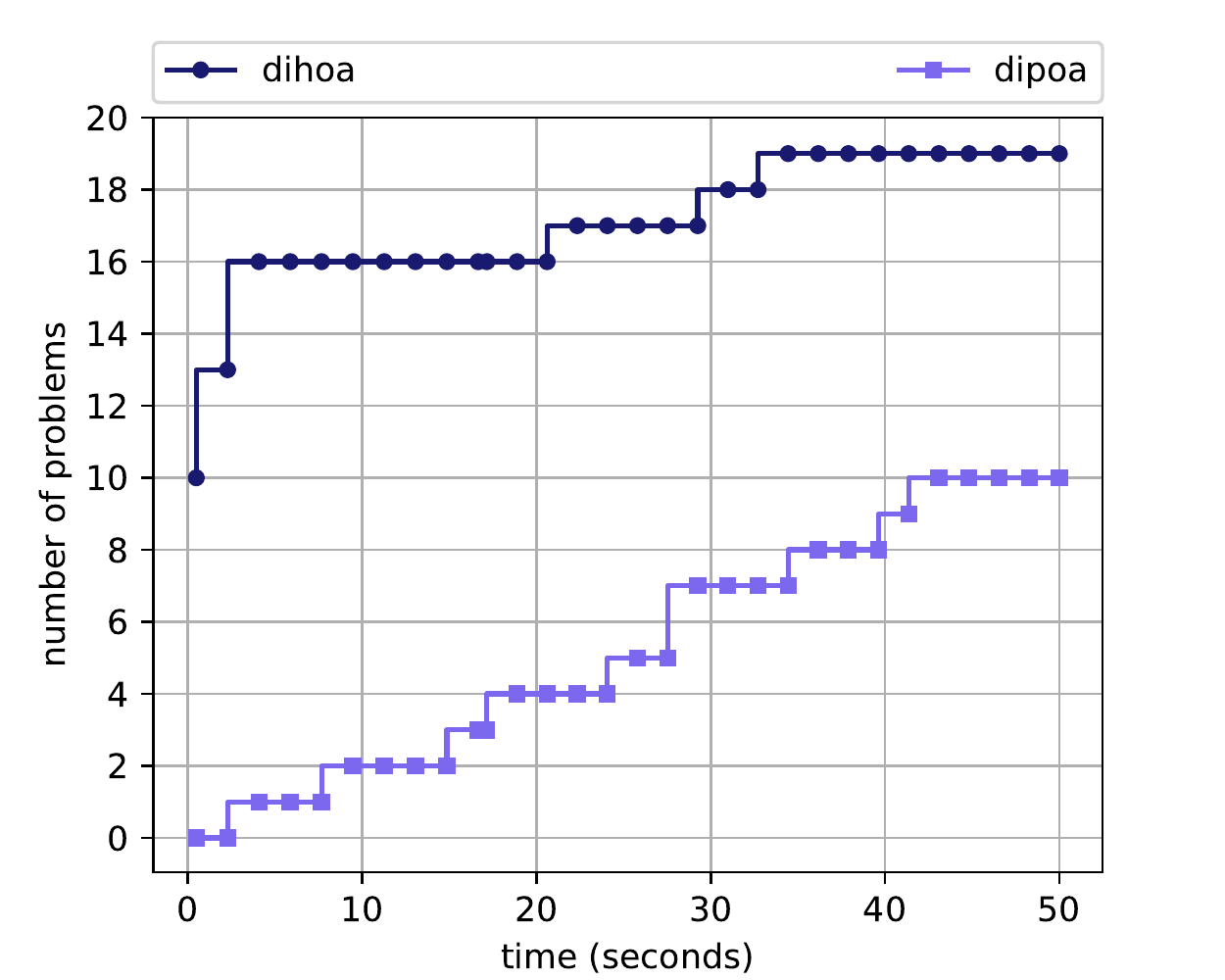}\label{fig:bench72}}}	
    \caption{Benchmark results for scenario $7$.}  \label{fig:bench7}
\end{figure}

This section presents performance profiles comparing the \vb{SCOT} performance with some MINLP solvers with different settings. We considered seven benchmark scenarios containing $20$ different problem instances with different properties and settings. In each scenario, we generate $15$ DSLogR and $5$ DSLinR random problem instances with a different number of features and sample points within a given range. Therefore the benchmark set consists of a total of $140$ problem instances. Moreover, each algorithm appearing in Table \ref{table:settings} is applied to solve all problem instances with $30$ different maximum execution times limits, starting from $0.5$ to $50$ seconds. Therefore, the total number of algorithm runs for each scenario is $600$, leading to $4200$ algorithm executions for all scenarios. Table \ref{table:scenarios} represents settings for each benchmark scenario where 
$n_{\min}$ and $n_{\max}$ are the minimum and the maximum number of features, $p_{\min}$ and $p_{\max}$ are the minimum and the maximum number of data points for each computational node, $n_p$ is the total number of problems, and $p_{tot}$ is the total number of data points considering all computational nodes.

Figures \ref{fig:bench1}-\ref{fig:bench7} compare the solvers \vb{SCOT}, \vb{Bonmin},  \vb{SHOT}, and \vb{Knitro} with both single-tree and multiple-tree algorithms. The comparison results for each scenario are shown as performance profiles consisting of two different sparsity levels of the solution. 

The benchmark results of the first three scenarios are depicted in Figures   \ref{fig:bench1}-\ref{fig:bench3} where \textit{small to medium size} problem instances are considered. In both sparsity levels, the DiPOA and DiHOA algorithms show better performance compared to other MINLP solvers. It can be observed in Figure \ref{fig:bench3} that for larger problem instances the performance gap between \vb{SCOT} and other MINLP solvers is increased. 

Figures \ref{fig:bench4}-\ref{fig:bench7} depict the performance profiles for the scenarios $4$-$7$ in which \textit{medium to large} problem instances are taken into account. In these scenarios, all MINLP solvers failed to solve the problem with the considered maximum execution time limit. Therefore, we only provide performance profiles of DiPOA and DiHOA algorithms. According to the performance profiles, the DiHOA algorithm is more efficient compared to DiPOA as the performance gap between them is increasing for large problem instances. For example, in scenario $7$ with $80\%$ of sparsity, DiHOA was able to solve $19$ problem instances within the given maximum execution time limit, whereas DiPOA only solved $10$ instances.
\begin{table}
\tbl{Benchmark Settings for Performance Profiles}
{\begin{tabular}{cccccccc} \toprule
 & \multicolumn{2}{l}{Scenario settings} \\ \cmidrule{2-8}
  {scenario} &{ $n_{\min}$} & { $n_{\max}$} & { $p_{\min}$} & { $p_{\max}$}&{ $N$} & { $n_p$} & { $p_{tot}$}  \\ \midrule
 {$1$}  & {$20$} & {$30$} &{$1000$} & {$2500$}& {$2$} &{$20$} & {$5000$} \\
 {$2$}  & {$20$} & {$30$} &{$1000$} & {$5000$}& {$2$} &{$20$} & {$10000$}\\
 {$3$}  & {$25$} & {$50$} &{$1000$} & {$5000$}& {$2$} &{$20$} & {$10000$}\\
 {$4$}  & {$25$} & {$100$} &{$1000$} & {$10000$}& {$4$} &{$20$} & {$40000$} \\
 {$5$}  & {$25$} & {$100$} &{$1000$} & {$20000$}& {$4$} &{$20$} & {$80000$} \\
 {$6$}  & {$25$} & {$100$} &{$1000$} & {$50000$}& {$4$} &{$20$} & {$200000$} \\ 
 {$7$}  & {$25$} & {$200$} &{$1000$} & {$50000$}& {$6$} &{$20$} & {$300000$}\\ \bottomrule
\end{tabular}}
\label{table:scenarios}
\end{table}

\subsection{Discussion}
In this section, we compared the performance of \vb{SCOT} using both DiHOA and DiPOA algorithms with some state-of-the-art MINLP solvers. According to the numerical experiments and performance profiles, the DiHOA algorithm achieves better performance and efficiency in all problem instances. However, one should keep in mind that SCOT is tailored for SCO problems which is not the case for the general purpose solvers. The results also showed that for large and distributed problems a distributed solver can provide an efficient and reliable solution.


\section{Conclusion}\label{sec:conclusion}
In this work, we introduced SCOT solver and DiHOA algorithm that SCOT implements to solve the DSCO problem \eqref{dccp}. The DiHOA implements an outer approximation algorithm to solve an MINLP equivalent to the SCO problem, combining single- and multiple-tree outer approximation with a decentralized algorithm for solving convex nonlinear subproblems, which generates primal solutions and cuts for the master dual problem.
The numerical benchmark results and comparison to state-of-the-art MINLP solvers indicated that SCOT equipped with DiHOA can be efficiently used for solving sparse convex optimization problems in different application domains. The performance and efficiency of the solver were achieved by augmenting SCOT with modern convex and mixed-integer optimization techniques. Continuing work aims at enhancing SCOT effectiveness by means of the development of decentralized mixed integer algorithms and new heuristic techniques to obtain outer approximations.


\section*{Disclosure statement}

The authors report there are no competing interests to declare.

\section*{Funding}
This work was funded in part by  Funda\c{c}\~ao de Amparo \`a Pesquisa e Inova\c{c}\~ao do Estado de Santa Catarina  (FAPESC) under grant 2021TR2265 and Coordena{\c{c}}{\~a}o de Aperfei{\c{c}}oamento de Pessoal de  N{\'i}vel Superior (CAPES, Brazil) under the project PrInt CAPES-UFSC 698503P1. Financial support from Digital Futures at KTH is also gratefully acknowledged.



\bibliographystyle{tfs} 
\bibliography{refs.bib}
\end{document}